\documentstyle[12pt,amssymb,amscd]{amsart}  

\newtheorem{thm}{Theorem}[subsection]

\newtheorem{cor}[thm]{Corollary}
\newtheorem{prop}[thm]{Proposition}
\newtheorem{lemma}[thm]{Lemma}

\newtheorem{conj}[thm]{Conjecture}

\theoremstyle{remark}
\newtheorem{remark}[thm]{Remark}

\theoremstyle{definition}

\numberwithin{equation}{section}

\newcommand{\bbA}{{\Bbb A}}
\newcommand{\bbC}{{\Bbb C}}
\newcommand{\bbH}{{\Bbb H}}
\newcommand{\bbL}{{\Bbb L}}
\newcommand{\bbR}{{\Bbb R}}
\newcommand{\bbW}{{\Bbb W}}
\newcommand{\bbZ}{{\Bbb Z}}
\newcommand{\cF}{{\cal F}}

\newcommand{\cO}{{\cal O}}
\newcommand{\cP}{{\cal P}}
\newcommand{\cL}{{\cal L}}
\newcommand{\cM}{{\cal M}}
\newcommand{\cN}{{\cal N}}
\newcommand{\cD}{{\cal D}}

\newcommand{\cA}{{\cal A}}

\newcommand{\cI}{{\cal I}}
\newcommand{\cC}{{\cal C}}
\newcommand{\cE}{{\cal E}}

\newcommand{\cR}{{\cal R}}
\newcommand{\cS}{{\cal S}}
\newcommand{\cT}{{\cal T}}

\newcommand{\chr}{\operatorname{char}}
\newcommand{\supp}{\operatorname{Supp}}
\newcommand{\SS}{\operatorname{SS}}

\newcommand{\Diff}{\operatorname{Diff}}
\newcommand{\Eu}{\operatorname{Eu}}
\newcommand{\eu}{\operatorname{eu}}
\newcommand{\D}{\operatorname{D}}
\newcommand{\R}{\operatorname{\bold R}}
\newcommand{\Hom}{\operatorname{Hom}}
\newcommand{\isomo}{\overset{\sim}{=}}
\newcommand{\isomoto}{\overset{\sim}{\to}}
\newcommand{\id}{\operatorname{id}}

\newcommand{\shHom}{\underline{\operatorname{Hom}}}
\newcommand{\colim}{\operatorname{colim}}
\newcommand{\holim}{\operatorname{holim}}
\newcommand{\Symm}{\operatorname{Sym}}

\begin{document}

\title{Riemann-Roch theorems via deformation quantization}

\author{P.Bressler}
\address{Department of Mathematics, The Pennsylvania State University,
University Park, PA 16802, USA}
\email{bressler@@math.psu.edu}
\author{R.Nest}
\address{Mathematical Institute, Universitatsparken, 5,
2100 Copenhagen, Denmark}
\email{rnest@@math.ku.dk}
\author{B.Tsygan}
\address{Department of Mathematics, The Pennsylvania State University,
University Park, PA 16802, USA}
\email{tsygan@@math.psu.edu}

\maketitle

\section{Introduction}
In this note we outline the proof of a formula expressing the Euler
class of a perfect complex of modules over a symplectic deformation
quantization of a complex manifold in terms of the Chern character of the 
associated symbol complex, the $\hat A$-class of the manifold and
a characteristic class of the deformation quantization. As a consequence we 
obtain the conjecture of P.Schapira and J.-P.Schneiders (\cite{SS},
Conjecture 8.5, see Conjecture \ref{conj:ch}) and the corollaries thereof.

In order to formulate the conjecture one needs a construction of
a local Chern character (cf. \cite{SS}, p.93). Such a construction
is provided in Section \ref{section:char} and appears to be new.

Specifically, we define
the Euler class and the Chern character of a perfect complex of sheaves
of modules over a sheaf of algebras by ``sheafifying'' the Dennis trace
map and the Goodwillie-Jones map defined by R.McCarthy in the generality
of exact categories in \cite{McC}. The consrtuction, as described in
Section \ref{section:char}, applies to sheaves of algebras on a topological
space, but, in fact, can be carried out for a sheaf of algebras on a site
(using, for example, the work of J.F.Jardine).

We reduce the Riemann-Roch theorem fot the Euler class to the Riemann-Roch
theorem for periodic cyclic cochains of deformed sheaves of algebras of
functions on a symplectic manifold. We reduce the latter to the local
Riemann-Roch theorem for periodic cyclic cochains of the Weyl algebra.
The latter result is contained in \cite{NT1}, \cite{NT2} and draws on the
ideas in \cite{FT2}.

For the sake of brevity we restrict our attention to the ``absolute'' case.
Analogous constructions can be carried out for families of symplectic
manifolds without any additional difficulties. Two techniques whose discussion
is not included here are
\begin{itemize}
\item Chern--Weil construction of Lie algebra cocycles with coefficients
``constant up to homotopy'' (used to reduce the Riemann-Roch theorem to
a local statement);
\item noncommutative differential calculus for periodic cyclic cochains
needed in the proof of the local Riemann-Roch theorem.
\end{itemize}
Admittedly, our presentation of the preliminaries is exteremely sketchy
A more detailed exposition will be given in a future paper. 

Here are the contents of the paper in brief. 

In Section \ref{section:ell} we recall some of the basic definitions and
facts about elliptic pairs leading up to Conjecture \ref{conj:ch}.

In Section \ref{section:quant} we review basic properties of quantized rings
of functions on symplectic manifolds and relate them to algebras of
differential and microdifferential operators. This relationship allows us
to place the traditional Riemann-Roch theorem into the context of deformation
quantization.

Section \ref{section:char} contains the background material
on Hochschild and cyclic homology, and their relationship to algebraic
K-theory necessary to define the Euler class and the Chern character of
a perfect complex of sheaves of modules over a sheaf of algebras.

In Section \ref{section:cx-man} we summarize the resuts on the Hochschild
and cyclic homology of algebras of (micro)differential operators as well as
quantized function algebras and compare characteristic classes and trace
density maps.

Section \ref{section:RR} is devoted to the statement of our main technical
result (Theorem \ref{thm:main}) which may be considered as a local
Riemann-Roch type theorem for periodic cyclic chains of the quantized ring
of functions. We apply Theorem \ref{thm:main} to characteristic classes of
perfect complexes of modules over the quantized rings of functions as well as
the algebras of differential and microdifferential operators to obtain
Conjecture \ref{conj:ch} and analogous statements as corollaries.

Section \ref{section:formalRR} leads up to the statement of Theorem
\ref{thm:formalRR} which is the analog of Theorem \ref{thm:main} in the
particular case when the symplectic manifold in question is the formal
neighborhood of the origin in a symplectic vector space.

In Section \ref{section:GF} we review the basic facts concerning Fedosov
connections and Gel'fand-Fuchs cohomology necessary to, loosely speaking,
establish a map from Theorem \ref{thm:formalRR} to Theorem
\ref{thm:main}. As a consequence of the existence of this (Gel'fand-Fuchs)
map the former theorem implies the latter.

The authors would like to thank J.-L.Brylinski and one another for
inspiring discussions.

\section{The index theory for elliptic pairs}\label{section:ell}
The application of our results to the index theory of elliptic pairs
is of particular importance. Here we recall the results of \cite{SS} restricting
ourselves to the ``absolute'' case for the sake of simplicity.

\subsection{}
Let $X$ be a complex manifold. An elliptic pair $(\cM^\bullet,F^\bullet)$
on $X$ consists
of an object $\cM^\bullet$ of $\D^b_{good}(\cD^{op}_X)$ (the bounded
derived category of complexes of right $\cD_X$-modules with {\em good}
\footnote{A $\cD_X$-module is called {\em good} if it admits a good filtration
in a neighborhood of every compact subset of $X$.}
cohomology and an object $F^\bullet$ of $\D^b_{\bbR-c}(X)$
(the bounded derived category of $\bbR$-constructible complexes of sheaves
of vector spaces over $\bbC$ on X) which satisfy
\[
\chr(\cM^\bullet)\cap\SS(F^\bullet)\subseteq T^*_XX\ .
\]

If
\[
\supp(\cM^\bullet,F^\bullet)\overset{def}{=}\supp\cM^\bullet\cap\supp F^\bullet
\]
(where the support of a complex of sheaves is understood to be the cohomological
support) is compact one has
\[
\dim H^\bullet(X;F^\bullet\otimes\cM^\bullet\otimes^{\bbL}_{\cD_X}\cO_X)
<\infty \ .
\]
Thus, the Euler characteristic
\[
\chi(X;F^\bullet\otimes\cM^\bullet\otimes^{\bbL}_{\cD_X}\cO_X)
\overset{def}{=}\sum_{i} (-1)^i
\dim H^i(X;F^\bullet\otimes\cM^\bullet\otimes^{\bbL}_{\cD_X}\cO_X)
\]
is defined. As particular cases one obtains
\begin{itemize}
\item the Euler characteristic $\chi(X;\cL^\bullet)$ of $X$ with coefficients
in a compactly supported perfect complex $\cL^\bullet$ of $\cO_X$-modules
(taking $F^\bullet = \bbC_X$, $\cM^\bullet = \cL^\bullet\otimes_{\cO_X}\cD_X$);

\item the index of an elliptic complex $\cL^\bullet$ on a compact real analytic
manifold $X_0$ (taking $X$ to be a complexification of $X_0$ such that
$\cL^\bullet$ extends to a complex $\tilde{\cL^\bullet}$ of $\cO_X$-modules and 
differential operators, $F^\bullet = \bbC_{X_0}$,
$\cM^\bullet = \Diff(\cO_X,\tilde{\cL^\bullet})$).
\end{itemize}

The index theorem for elliptic pairs (\cite{SS}, Theorem 5.1) says that,
for an elliptic pair $(\cM^\bullet,F^\bullet)$ with compact support on $X$
of dimension $\dim_{\bbC}X = d$
\[
\chi(X;F^\bullet\otimes\cM^\bullet\otimes^{\bbL}_{\cD_X}\cO_X)
=
\int_{T^*X}\mu\eu(\cM^\bullet)\smile\mu\eu(F^\bullet)
\]
where $\mu\eu(\cM^\bullet)\in H^{2d}_{\chr(\cM^\bullet)}(T^*X;\bbC)$ and 
$\mu\eu(F^\bullet)\in H^{2d}_{\SS(F^\bullet)}(T^*X;\bbC)$ are defined in
\cite{SS}.

The class $\mu\eu(F^\bullet)$ is the characteristic cycle of the constructible
complex $F^\bullet$ as defined by Kashiwara (see \cite{KS} for more details).
For example, if $Y\subset X$ is a closed real analytic submanifold, one has
$\mu\eu(\bbC_Y) = [T^*_YX]$.

With regard to $\mu\eu(\cM^\bullet)$ P.Schapira
and J.-P.Schneiders conjectured that it is related to a certain 
characteristic class of the symbol $\sigma(\cM^\bullet)$ of $\cM^\bullet$
(\cite{SS}, Conjecture 8.5, see Conjecture \ref{conj:ch}).

\subsection{The Riemann-Roch type formula}
Suppose that $\cM^\bullet$ is an object of $\D^b_{good}(\cD^{op}_X)$.
Let $\pi : T^*X\to X$ denote the projection.

If $\cM^\bullet$ admits a gobal {\em good} filtration (and this is the case
when $X$ is compact) then the symbol complex
of $\cM$ is defined by
\[
\sigma(\cM^\bullet) = \pi^{-1}gr\cM^\bullet\otimes_{\pi^{-1}gr\cD_X}
\cO_{T^*X}\ .
\]
The assumption that the filtration is good amounts to the fact that
$\sigma(\cM^\bullet)$ has $\cO_{T^*X}$-coherent cohomology, and the
characteristic variety is defined by
\[
\chr(\cM^\bullet)\overset{def}{=}\supp\sigma(\cM^\bullet)\ .
\]
For $\Lambda$ a closed subvariety of $T^*X$ let
$K^0_\Lambda(T^*X)$ denote the Grothendieck group of perfect complexes of
$\cO_{T^*X}$-modules supported on $\Lambda$ (i.e. acyclic on the complement
of $\Lambda$ in $T^*X$).

For $\Lambda$ containing $\chr(\cM)$ let 
$\sigma_\Lambda(\cM^\bullet)$ denote the class of $\sigma(\cM^\bullet)$
in $K^0_\Lambda(T^*X)$.

\begin{remark}
As is easy to show, both the characteristic variety
and the class of the symbol in the Grothendieck group are independent of
the choice of the good filtration, thus the existence of a good filtration
locally is sufficient to define $\sigma_\Lambda(\cM^\bullet)$.
\end{remark}

One can define the Chern character (see \eqref{map:ch})
\[
ch_\Lambda : K^0_\Lambda(T^*X)\to\bigoplus_i H^{2i}_\Lambda(T^*X;\bbC)\ .
\]

For $\alpha$ an element of a graded object let $[\alpha]^p$ denote the
homogeneous component of $\alpha$ of degree $p$.

In \cite{SS}, P.Schapira and J.-P.Schneiders make the following conjecture.

\begin{conj}\label{conj:ch}
For $\cM^\bullet$ in $\D^b_{good}(\cD^{op}_X)$, $\Lambda$ a conic subvariety
of $T^*X$ containing $\chr(\cM^\bullet)$
\[
\mu\eu(\cM^\bullet) =
\left[ ch_\Lambda(\sigma(\cM^\bullet)\smile \pi^*Td(TX)\right]^{2d}
\]
\end{conj}

We will obtain the above conjecture as a corollary of a Riemann-Roch type
formula in the context of deformation quantization of symplectic manifolds.

\section{Deformation quantization}
\label{section:quant}
\subsection{Review of deformation quantization}
A {\em deformation quantization} of a manifold $M$ is a formal one parameter
deformation of the structure sheaf $\cO_M$, i.e. a sheaf of algebras
$\bbA^\hbar_M$ flat over $\bbC [[\hbar]]$ together with an isomoprhism of
algebras $\bbA^\hbar_M\otimes_{\bbC [[\hbar]]}\bbC\to\cO_M$.

The formula
\[
\lbrace f, g\rbrace = \frac{1}{\hbar}[\tilde f,\tilde g ]
+ \hbar\cdot \bbA^\hbar_M\ ,
\]
where $f$ and $g$ are two local sections of $\cO_M$ and $\tilde f$, $\tilde g$
are their respective lifts $\bbA^\hbar_M$, defines a Poisson structure on $M$
called the Poisson structure associated to the deformation quantization
$\bbA^\hbar_M$.

\begin{remark}
The definition of ``deformation quantization'' as given above is essentially
the one given by \cite{BFFLS} in the case of $C^\infty$ manifolds. It is not
at all clear whether the scope of generality of the above definition above
is sufficiently broad ... 
\end{remark}

The deformation quantization $\bbA^\hbar_M$ is called {\em symplectic} if the 
associated Poisson structure is nondegenerate. In this case $M$ is symplectic.
In what follows we will only consider symplectic deformation quantizations,
so assume that $\bbA^\hbar_M$ is symplectic from now on.

It is known (and not difficult to show) that all symplectic deformation
quantizations of $M$ of dimension $\dim_{\bbC}M = 2d$ are locally isomorphic
to the standard deformation quantization of $\bbC^{2d}$. That is, for any
point $x\in M$ and small neighborhoods $U$ of $x$ and $U'$ of the origin in 
$\bbC^{2d}$ there is an isomorphism
\begin{equation}\label{loc-iso-h}
\bbA^\hbar_{\bbC^{2d}}(U')\overset{def}{=}
\cO_{\bbC^{2d}}(U')[[\hbar ]]\isomoto\bbA^\hbar_M(U)
\end{equation}
of algebras over $\bbC [[\hbar ]]$, continuous in the $\hbar$-adic topology,
where the product on $\bbA^\hbar_{\bbC^{2d}}(U')$ is given, in coordinates
$x_1,\ldots,x_d,\xi_1,\ldots,\xi_d$ on $\bbC^{2d}$ by the standard Weyl product
\[
(f\ast g)(\underline x,\underline\xi) = \\
exp\left(\frac{\sqrt{-1}\hbar}{2}\sum_{i=1}^d\left(
\frac{\partial\ }{\partial\xi_i}\frac{\partial\ }{\partial y_i}-
\frac{\partial~}{\partial\eta_i}\frac{\partial~}{\partial x_i}\right)\right)
f(\underline x,\underline\xi)g(\underline y,\underline\eta)
\vert_{\overset{\underline x=\underline y}{\underline\xi =\underline\eta}}
\]
where $\underline x = (x_1,\ldots,x_d),\ \underline\xi = (\xi_1,\ldots,\xi_d),
\ \underline y = (y_1,\ldots,y_d),\ \underline\eta = (\eta_1,\ldots,\eta_d)$.
Note that the reduction of \eqref{loc-iso-h} modulo $\hbar$ is an isomorphism
of Poisson algebras
\[
\cO_{\bbC^{2d}}(U')\isomoto\cO_M(U)
\]
with the Poisson brackets associated to the standard symplectic structure on
$\bbC^{2d}$ and to the deformation quantization respectively. In particular
(the images of) $x_1,\ldots,x_d,\xi_1,\ldots,\xi_d$ form a Darboux coordinate
system on $U\subset M$.

To a symplectic deformation quantization $\bbA^\hbar_M$ one associates
a characteristic class $\theta\in H^2(M;\frac{1}{\sqrt{-1}\hbar}\bbC[[\hbar]])$
with the property that the coefficient of $\frac{1}{\sqrt{-1}\hbar}$ is the
class of the symplectic form associated to the deformation quantization.

\subsection{Microlocalization}\label{subsec:microloc}
We start with a coherent $\cD_X$-module
$\cM$ equipped with a good filtration $F_\bullet\cM$ and consider
the graded module
\[
\cR\cM\overset{def}{=} \bigoplus_q F_q\cM\cdot\hbar^q
\]
over the graded ring
\[
\cR\cD_X \overset{def}{=} \bigoplus_q F_q\cD_X\cdot\hbar^q
\subset \cD_X [\hbar ]
\]
where $F_\bullet\cD_X$ is the filtration by the order of the differential
operator.

There exists a deformation quantization $\bbA^\hbar_{T^*X}$ of $T^*X$,
i.e. a formal deformation of the structure sheaf $\cO_{T^*X}$ and
faithfully flat maps
\[
\pi^{-1}\cR\cD_X @>>> \cR\cE_X @>>> \bbA^\hbar_{T^*X}
\]
of algebras over $\bbC [\hbar ]$ where $\cE_X$ is the sheaf (on $T^*X$) of microdifferential operators. The characteristic class $\theta$ of the 
deformation $\bbA^\hbar_{T^*X}$ is equal to
$\displaystyle\frac12\pi^*c_1(TX)$ (note that the sympledctic form is exact in
this case).

Consider the ``microlocalization''
\[
\mu\cM\overset{def}{=}\pi^{-1}\cR\cM\otimes_{\pi^{-1}\cR\cD_X}\bbA^\hbar_{T^*X}
\]
of the filtered $\cD_X$-module $\cM$. Note that $\mu\cM$ is $\hbar$-torsion
free. Then, clearly, there is an isomorphism
\[
\sigma(\cM)\isomo\mu\cM\otimes_{\bbA^\hbar_{T^*X}}\cO_{T^*X}\isomo
\mu\cM/\mu\cM\cdot\hbar\ .
\]

Defining the symbol of an $\hbar$-torsion free $\bbA^\hbar_{T^*X}$-module
$\cN$ by
\[
\sigma(\cN) = \cN\otimes_{\bbA^\hbar_{T^*X}}\cO_{T^*X}
\]
we have
\[
\sigma(\cM)\isomo\sigma(\mu\cM)\ .
\]

\section{Characteristic classes of perfect complexes}\label{section:char}
In this section we construct the Euler class (with values in Hochschild
homology) and the Chern character (with values in negative cyclic homology)
for a perfect complex of sheaves of modules over a sheaf of algebras over
a topological space.

\subsection{Review of Hochschild and cyclic homology}
Let $k$ denote a commutative algebra over a field of characteristic zero
and let $A$ be a flat $k$-algebra with $1_A\cdot k$ contained in the center,
not necessarily commutative. Let
$C_p(A)\overset{def}{=}A^{\otimes_k p+1}$ and let
\begin{eqnarray*}
b : C_p(A) & @>>> & C_{p-1}(A) \\
a_0\otimes\cdots\otimes a_p & \mapsto & 
(-1)^p a_pa_0\otimes\cdots\otimes a_{p-1} + \\ & &
\sum_{i=0}^{p-1}(-1)^ia_0\otimes\cdots\otimes a_ia_{i+1}\otimes\cdots\otimes
a_p\ .
\end{eqnarray*}
Then $b^2 = 0$ and the complex $(C_\bullet, b)$, called {\em the standard
Hochschild complex of $A$} represents $A\otimes^{\bbL}_{A\otimes_kA^{op}}A$
in the derived category of $k$-modules.

The map
\begin{eqnarray*}
B : C_p(A) & @>>> & C_{p+1}(A) \\
a_0\otimes\cdots\otimes a_p & \mapsto & \sum_{i=0}^p (-1)^{pi}
1\otimes a_i\otimes\cdots\otimes a_p\otimes a_0\otimes\cdots\otimes a_{i-1}
\end{eqnarray*}
satisfies $B^2 = 0$ and $[B,b] = 0$ and therefore defines a map of complexes
\[
B: C_\bullet(A) @>>> C_\bullet(A)[-1]\ .
\]
For $i,j,p\in\bbZ$ let
\begin{eqnarray*}
CC^-_p(A) & = & \prod_{\overset{i\geq 0}{i+j = p \mod 2}} C_{i+j}(A) \\
CC^{per}_p(A) & = & \prod_{i+j = p \mod 2} C_{i+j}(A)\ .
\end{eqnarray*}
The complex $(CC^-_\bullet(A),B+b)$ (respectively
$(CC^{per}_\bullet(A),B+b)$) is called the {\em negative} (respectively
{\em periodic}) {\em cyclic complex of $A$}.

There are inclusions of complexes
\[
CC^-_\bullet(A)[-2]\hookrightarrow CC^-_\bullet(A)\hookrightarrow
CC^{per}_\bullet(A)
\]
and the short exact sequence
\[
0 @>>> CC^-_\bullet(A)[-2] @>>> CC^-_\bullet(A) @>>> C_\bullet(A) @>>> 0\ .
\]

Suppose that $X$ is a topological space and $\cA$ is a flat sheaf of 
$k$-algebras on $X$ such that there is a global section $1\in\Gamma(X;\cA)$ 
which
restricts to $1_{\cA_x}$ and $1_{\cA_x}\cdot k$ is contained in the center of 
$\cA_x$ for every point $x\in X$. Let $C_\bullet(\cA)$ (respectively
$CC^-_\bullet(\cA),\ CC^{per}_\bullet(\cA)$) denote the complex
of sheaves of $k$-modules associated to the presheaf with value
$C_\bullet(\cA(U))$ (respectively
$CC^-_\bullet(\cA(U)),\ CC^{per}_\bullet(\cA(U))$) on an open subset $U$ of
$X$. Then $C_\bullet(\cA)$ represents
$\cA\otimes^{\bbL}_{\cA\otimes_k\cA^{op}}\cA$ in the derived category of
sheaves of $k$-modules on $X$.

\subsection{Perfect complexes}
We briefly recall the notion of perfection as introduced in \cite{Ill}.

A module over $\cA$ is said to be {\em free of finite type} if it is
isomorphic to $\cA^{\oplus n}$ for some $n\in\bbZ$.

A complex $\cP^\bullet$ of (sheaves of) $\cA$-modules is called {\em strictly
perfect} if
\begin{enumerate}
\item $\cP^p = 0$ for almost all $p\in\bbZ$;

\item for any $p\in\bbZ$ and every point $x\in X$ there exists a neighborhood
$U$ of $x$ such that $\cP^p\vert_U$ is a direct summand of a free
$\cA\vert_U$-module of finite type.
\end{enumerate}

A complex $\cF^\bullet$ of (sheaves of) $\cA$-modules is called {\em perfect}
if for any point $x\in X$ there exists an open neighborhood $U$ of $x$,
a strictily perfect complex $\cP^\bullet$ of $\cA\vert_U$-modules and
a quasiisomorphism $\cP^\bullet\to\cF^\bullet\vert_U$.

\subsection{The Euler class in Hochschild homology}
For a perfect
complex $\cF^\bullet$ of $\cA$ modules we define the {\em Lefschetz map}
as the morphism in the derived category of sheaves of $k$-modules given
by the composition
\[
\R\shHom^\bullet_{\cA}(\cF^\bullet,\cF^\bullet)
@<{\isomo}<< (\R\shHom^\bullet_{\cA}(\cF^\bullet,\cA)\otimes_k\cF^\bullet)
\otimes^{\bbL}_{\cA\otimes_k\cA^{op}}\cA @>{ev\otimes\id}>>
\cA\otimes^{\bbL}_{\cA\otimes_k\cA^{op}}\cA
\]
and will denote it by $\cL_{\cA}(\cF^\bullet)$ or simply by $\cL_{\cA}$.

We define the {\em Euler map} of $\cF^\bullet$ as the morphism in
the derived category given by the composition
\[
k @>{1\mapsto\id}>> \R\shHom^\bullet_{\cA}(\cF^\bullet,\cF^\bullet)
@>{\cL_{\cA}}>> \cA\otimes^{\bbL}_{\cA\otimes_k\cA^{op}}\cA
\]
and will denote it by $\Eu_{\cA}(\cF^\bullet)$.

Suppose that $Z$ is a closed subset of $X$ such that
$Z\supseteq\supp\cF^\bullet\overset{def}{=}\bigcup_p\supp H^p\cF^\bullet$.
Then, clearly, the canonical morphism 
$\R\Gamma_Z(\R\shHom^\bullet_{\cA}(\cF^\bullet,\cF^\bullet))\to
\R\shHom^\bullet_{\cA}(\cF^\bullet,\cF^\bullet)$ is an isomorphism.

Thus, after applying the functor $\R\Gamma(X;\R\Gamma_Z(\bullet))$
and passing to
cohomology, the morphism $\Eu_{\cA}(\cF^\bullet)$ determines a cohomology
class in $H^0_Z(X;\cA\otimes^{\bbL}_{\cA\otimes_k\cA^{op}}\cA)$ which we 
will refer to as the {\em Euler class} and will denote by
$\eu_{\cA}^Z(\cF^\bullet)$.

For $Z$ a closed subset of $X$ let $K^0_Z(\cA)$ denote the Grothendieck group
of perfect complexes of $\cA$-modules supported on $Z$; set $K^0(\cA)
\overset{def}{=}K^0_X(\cA)$. The Euler class defined above determines
a homomorphism of groups
\[
\eu_{\cA}^Z : K^0_Z(\cA) @>>>
H^0_Z(X;\cA\otimes^{\bbL}_{\cA\otimes_k\cA^{op}}\cA)\ .
\]

\subsection{Review of K-theory}
Recall that, to a category $\cC$ with cofibrations and weak equivalences
one associates the simplicial category $S_\bullet\cC$ (the $S$-construction)
whose definition can be found in \cite{W}.
The $K$-theory spectrum $K(\cC)$ of $\cC$ is defined by
$K(\cC) = \Omega\vert S_\bullet\cC\vert$ which describes the zeroth
$\Omega^\infty$-space whose deloopings are given by iterating the
$S$-construction. The groups $K_i(\cC)$ are defined by
\[
K_i(\cC) = \pi_iK(\cC)\ .
\]

To a ring $A$ one associates the spectrum $K^{naive}(A)$ (respectively $K(A)$)
defined as the $K$-theory spectrum of the category of strictly perfect
(respectively perfect) complexes of $A$-modules.

The natural inclusion of the former category into the latter induces the
natural morphism of spectra
\[
K^{naive}(A) @>>> K(A)\ .
\]

The inclusion of the category $\cP_A$ of finitely generated projective
$A$-modules into the category of strictly perfect complexes induces
the natural weak equivalence
\[
K(\cP_A) @>>> K^{naive}(A)\ .
\]

\subsection{Hochschild and cyclic homology revisited}
We summarize some of the results of \cite{McC} concerning the Hochschild
and cyclic homology of exact categories.

To a $k$-linear additive category $\cC$ ($k$ a ring) one associates {\em
the cyclic nerve} $CN_\bullet\cC$ which is a cyclic $k$-module.

In the case when $\cC$ is a $k$-linear additive category with one object,
i.e. a $k$-algebra $A$, one has $CN_\bullet\cC\isomo C_\bullet(A)$.

In what follows we will not make notational distinctions between complexes
and corresponding simplicial Abelian groups.

To a category $\cC$ with cofibrations and weak equivalences one associates
the Hochschild complex $C_\bullet(\cC)$ by
\[
C_\bullet(\cC) = \Omega\vert CN_\bullet S_\bullet\cC\vert =
\operatorname{Tot}(CN_\bullet S_\bullet\cC)[-1]\ .
\]
Using the cyclic structure on $C_\bullet(\cC)$ one defines the cyclic
(respectively the negative cyclic, respectively the periodic cyclic) complex
of $\cC$ which we denote by $CC_\bullet(\cC)$ (respectively
$CC^-_\bullet(\cC)$, respectively $CC^{per}_\bullet(\cC)$). There is
natural commutative diagram
\[
\begin{CD}
CC^-_\bullet(\cC) @>>> CC^{per}_\bullet(\cC) @>>> CC_\bullet(\cC)[2] \\
@VVV			@VVV				@VV{\id}V  \\
C_\bullet(\cC) @>>> CC_\bullet(\cC) @>>> CC_\bullet(\cC)[2]
\end{CD}
\]
with rows exact triangles.

Let $\cP_A$ denote the category of finitely generated projective
$A$-modules. There is a natural quasiisomorphism of cyclic $k$-modules
\[
C_\bullet(A) @>>> C_\bullet(\cP_A)
\]
which induces quasiisomorphisms of respective cyclic (negative cyclic, periodic
cyclic) complexes.

The natural inclusion
\[
S_\bullet\cC @>>> CN_\bullet S_\bullet\cC
\]
induces the natural morphism of spectra
\[
\Eu_{\cC} : K(\cC) @>>> C_\bullet(\cC)
\]
which is the Dennis trace map in the case $\cC = \cP_A$.
The morphism $\Eu$ has a natural lifting
\[
ch_{\cC} : K(\cC) @>>> CC^-_\bullet(\cC)
\]
which is the Chern character (the Goodwillie-Jones map) in the case
$\cC = \cP_A$.

One can show that the (exact) inclusion of $\cP_A$ into the category of
strictly perfect complexes of $A$-modules induces quasiisomorphisms
of respective Hochschild, cyclic, negative cyclic, and periodic cyclic
complexes. In particular one has the Euler class and the Chern character
\[
\Eu_A : K^{naive}(A) @>>> C_\bullet(A),\ \ \ 
ch_A : K^{naive}(A) @>>> CC^-_\bullet(A)\ .
\]

\subsection{Presheaves of spectra}
The following is a summary of the definitions and constructions of Section 3
of \cite{M}. We refer the reader to \cite{M} (and references therein,
particularly \cite{T}) for further details.

To a presheaf $\cS$ of spectra on $X$ (and, more generally, a presheaf with
values in a category with filtered colimits, products and coproducts) one can
associate the functorial cosimplicial Godement resolution in the usual way.
The stalk $\cS_x$ of $\cS$ at $x\in X$ is defined by
\[
\cS_x = \colim \cS(U)
\]
where the colimit is taken over all open neighborhoods $U$ of $x$.
For $U$ an open subset of $X$ let
\[
T(\cS)(U) = \prod_{x\in U}\cS_x\ .
\]
Then $T(\cS)$ is a presheaf on $X$, $T$ is a functor from presheaves to
presheaves and in fact a monad. Thus, $T$ gives rise to the functorial
(in $\cS$) cosimplicial object $T^\bullet(\cS)$ in presheaves on $X$
called the {\em Godement resolution of} $\cS$.

For $\cS$ a presheaf of spectra on $X$ and an open subset $U$ of $X$ let
\[
\bbH(U;\cS) = \underset{\Delta}{\holim}\ T^\bullet(\cS)(U)\ .
\]
The assignment $U\mapsto\bbH(U;\cS)$ determines a presheaf of spectra on $X$.
Note that there is a natural morphism of presheaves of spectra
\[
\cS @>>> \bbH(\bullet ;\cS)\ .
\]

For $Z\subset X$ a closed subset let $\Gamma_Z\cS$ denote the presheaf of
spectra whose value on an open set $U\subset X$ is defined to be the homotopy
fiber of the restriction map $\cS(U)\to \cS(U\setminus Z)$. In particular
there is a canonical morphism $\Gamma_Z\cS\to\cS$. Set
\[
\bbH_Z(U;\cS)\overset{def}{=}
\underset{\Delta}{\holim}\ \Gamma_ZT^\bullet(\cS)(U)\ .
\]

\subsection{The Euler class and the Chern character revisited}
Let $K_Z(\cA)$ (respectively $K^{naive}_Z(\cA)$) denote the $K$-theory
spectrum of the category of perfect complexes (respectively strictly perfect
complexes) of sheaves of $\cA$-modules on $X$ supported on $Z$ (i.e. acyclic
on the complement of $Z$ in $X$). The inclusion of the category of strictly 
perfect complexes into the category of all perfect complexes induces the
natural morphism of spectra $K^{naive}_Z(\cA)\to K_Z(\cA)$.

The assignment $U\mapsto K^{naive}(\cA(U))$ (respectively
$U\mapsto K^{naive}(\cA\vert_U)$, $U\mapsto K(\cA\vert_U)$) determines a
presheaf of spectra on $X$. The functor $M\mapsto M\otimes_{\cA(U)}\cA\vert_U$
and the inclusion of strictly perfect complexes into perfect complexes induce
morphism of presheaves of spectra
\begin{equation}\label{maps:KKK}
K^{naive}(\cA(\bullet)) @>>> K^{naive}(\cA\vert_\bullet)
@>>> K(\cA\vert_\bullet)\ .
\end{equation}

\begin{lemma}
The morphisms \eqref{maps:KKK} induce equivalences on stalks.
\end{lemma}

\begin{cor}
The morphism \eqref{maps:KKK} induce equivalences
\begin{equation}\label{maps:HHH}
\bbH_Z(X;K^{naive}(\cA(\bullet))) @>>>
\bbH_Z(X;K^{naive}(\cA\vert_\bullet)) @>>>
\bbH_Z(X;K(\cA\vert_\bullet))\ .
\end{equation}
\end{cor}

Combing the equivalences \eqref{maps:HHH} with the the Euler class
(Dennis trace map) and the canonical morphism
$K_Z(\cA)\to\bbH_Z(X;K(\cA\vert_\bullet))$ we obtain the morphism of spectra
\[
\Eu :K_Z(\cA) @>>> \bbH_Z(X;C_\bullet(\cA)) = \R\Gamma_Z(X;C_\bullet(\cA))\ .
\] 
Using the Chern character instead of the Euler class we obtain the morphism
\begin{equation}\label{map:ch}
ch :K_Z(\cA) @>>> \bbH_Z(X;CC^-_\bullet(\cA)) =
\R\Gamma_Z(X;CC^-_\bullet(\cA))\ .
\end{equation} 

All in all, we have constructed Euler class and Chern character for a perfect
complex.

\section{Characteristic classes on complex manifolds}
\label{section:cx-man}
\subsection{Notations and conventions}
In what follows we will be considering, for a complex
manifold $X$, the sheaves of algebras
$\cO_X$, $\cD_X$, $\cE_{T^*X}$, and $\bbA^\hbar_{T^*X}$. All of these are
sheaves of topological vector spaces. In what follows all tensor products
are understood to be projective tensor products. In particular, let
\begin{eqnarray*}
\cO_X^{\frak{e}} & = & \cO_X\widehat\otimes\cO_X \\
\cD_X^{\frak{e}} & = & \cD_X\widehat\otimes\cD_X^{op} \\
\cE_X^{\frak{e}} & = & \cE_X\widehat\otimes\cE_X \\
(\bbA^\hbar_M)^{\frak{e}} & = & \bbA^\hbar_M
\widehat\otimes_{\bbC[[\hbar]]}(\bbA^\hbar_M)^{op}\ .
\end{eqnarray*}

If $\cA$ is one of $\cO_X$, $\cD_X$, $\cE_X$, $\bbA^\hbar_M$ then
$\cA$ has a natural structure of an $\cA^{\frak{e}}$-module and there
is a natural map $\cA\otimes\cA^{op}\to\cA^{\frak{e}}$ which induces a map
of complexes
\[
\cA\otimes^{\bbL}_{\cA\otimes_k\cA^{op}}\cA @>>>
\cA\otimes^{\bbL}_{\cA^{\frak{e}}}\cA\ .
\]
In what follows we will only consider the composition
\[
\R\shHom^\bullet_{\cA}(\cF^\bullet,\cF^\bullet) @>{\cL_{\cA}}>>
\cA\otimes^{\bbL}_{\cA\otimes_k\cA^{op}}\cA @>>>
\cA\otimes^{\bbL}_{\cA^{\frak{e}}}\cA\ ,
\]
refer to it as the Lefschetz map and denote it by $\cL_{\cA}$. Similarly,
the Euler map $\Eu_{\cA}(\cF^\bullet)$ will denote the composition
\[
k @>{1\mapsto\id}>> \R\shHom^\bullet_{\cA}(\cF^\bullet,\cF^\bullet)
@>{\cL_{\cA}}>> \cA\otimes^{\bbL}_{\cA^{\frak{e}}}\cA\ .
\]

\subsection{Some examples of Hochschild and cyclic homology}
Here we recall some facts about the Hochschild and cyclic homology of
$\cO_X$, $\cD_X$, $\cE_X$ and $\bbA^\hbar_{T^*X}$ for a complex manifold $X$
of dimension $\dim_{\bbC}X = d$ from \cite{B}.

By the theorem of Hochschild-Kostant-Rosenberg there is an isomorphism
in the derived category of $\cO_X$-modules
\[
\mu_{\cO_X} : \cO_X\otimes^{\bbL}_{\cO_X^{\frak{e}}}\cO_X
@>>> \bigoplus_p\Omega^p_X [p]
\]
given, in terms of the standard Hochschild complex representing
$\cO_X\otimes^{\bbL}_{\cO_X\otimes\cO_X}\cO_X$, by the formula
\[
a_0\otimes\cdots\otimes a_p \mapsto \frac{1}{p!}
a_0 da_1\wedge\cdots\wedge da_p \ .
\]
The same formula gives the quaiisomorphism
\begin{equation}\label{map:HKR}
\tilde\mu_{\cO} : CC^{per}_\bullet(\cO_X) @>>>
\prod_{p\in\bbZ}\Omega^\bullet_X [2p]
\end{equation}
of the periodic cyclic complex of $\cO_X$ and the de Rham complex of $X$,
since it is easily verified that the diagram
\[
\begin{CD}
C_\bullet(\cO_X) @>{B}>> C_\bullet(\cO_X)[-1] \\
@V{\mu}VV			@V{\mu}VV \\
\bigoplus_p\Omega^p_X [p] @>{d}>> \bigoplus_p\Omega^p_X [p-1]
\end{CD}
\]
(where $d$ denotes exterior differentiation) is commutative.

In view of the fact that the canonical map $\bbC_X\to\Omega^\bullet_X$
is a quasiisomorphism we will view the map $\tilde\mu_{\cO}$ as the isomorphism
(in the derived category)
\[
\tilde\mu_{\cO} : CC^{per}_\bullet(\cO_X) @>>>
\prod_{p\in\bbZ}\bbC_X [2p]\ .
\]
The inverse to $\tilde\mu_{\cO}$ is provided by the map of periodic cyclic
complexes induced by the inclusion $\bbC_X\hookrightarrow\cO_X$.

According to J.-L.Brylinski (\cite{Bry})
\[
H^p(\cD_X\otimes^{\bbL}_{\cD_X^{\frak{e}}}\cD_X)\isomo
\left\lbrace\begin{array}{lll} 0 & \text{if} & p\neq -2d \\
\bbC_X & \text{if} & p = -2d \end{array}\right.
\]
Consider a point $x\in X$, an open neighborhood $U$ of $x$ in $X$ and a local 
coordinate system $x_1,\ldots,x_d$ centered at $x$. The class of the Hochschild 
cycle $Alt(1\otimes x_1\otimes\cdots x_d\otimes
\displaystyle\frac{\partial}{\partial x_1}\otimes\cdots\otimes
\displaystyle\frac{\partial}{\partial x_d})$ represents a global section
$\Phi^{\cD}_U$ of $H^{-2d}(\cD_U\otimes^{\bbL}_{\cD_U^{\frak{e}}}\cD_U)$. Let
\[
\mu_{\cD} : \cD_X\otimes^{\bbL}_{\cD_X^{\frak{e}}}\cD_X
\to \bbC_X [2d]
\]
denote the isomorphism in the derived category of sheaves which corresponds
to the global section of $H^{-2d}(\cD_X\otimes^{\bbL}_{\cD_X^{\frak{e}}}\cD_X)$
determined by the condition that it restricts to $\Phi^{\cD}_U$ for every 
sufficiently small open set $U\subset X$. Let
\[
\mu^\hbar_{\cD} : \cD_X\otimes^{\bbL}_{\cD_X^{\frak{e}}}\cD_X
\to \bbC_X [\hbar^{-1},\hbar]] [2d]
\]
denote the composition $\displaystyle\frac{1}{\hbar^d}\mu_{\cD}$.

Similar result hold for the sheaf $\cE_X$ of microdifferential operators.
Specifically,
\[
H^p(\cE_X\otimes^{\bbL}_{\cE_X^{\frak{e}}}\cE_X)\isomo
\left\lbrace\begin{array}{lll} 0 & \text{if} & p\neq -2d \\
\bbC_{T^*X} & \text{if} & p = -2d \end{array}\right.
\]
Let $\xi_i$ denote the symbol of $\displaystyle\frac{\partial}{\partial x_i}$
where
$x_1,\ldots,x_d$ are local coordinates on $X$ as before. The class of
the Hochschild cycle
$Alt(1\otimes x_1\otimes\cdots\otimes x_d\otimes\xi_i\otimes\cdots\otimes\xi_d)$
determines a global
section $\Phi^{\cE}_U$ of 
$H^{-2d}(\cE_U\otimes^{\bbL}_{\cE_U^{\frak{e}}}\cE_U)$. Let
\[
\mu_{\cE} : \cE_X\otimes^{\bbL}_{\cE_X^{\frak{e}}}\cE_X
\to \bbC_{T*X} [2d]
\]
denote the isomorphism in the derived category of sheaves which corresponds
to the global section of $H^{-2d}(\cE_X\otimes^{\bbL}_{\cE_X^{\frak{e}}}\cE_X)$
determined by the condition that it restricts to $\Phi^{\cE}_U$ for every 
sufficiently small open set $U\subset X$. Let
\[
\mu^\hbar_{\cE} : \cE_X\otimes^{\bbL}_{\cE_X^{\frak{e}}}\cE_X
\to \bbC_{T^*X} [\hbar^{-1},\hbar]] [2d]
\]
denote the composition $\displaystyle\frac{1}{\hbar^d}\mu_{\cE}$.

Suppose that $\bbA^\hbar_M$ is a symplectic deformation quantization of
a complex manifold $M$ of dimension $\dim_{\bbC} M = 2d$. The sheaf of
algebras $\bbA^\hbar_M[\hbar^{-1}]$ exhibits properties similar to those
of $\cE_X$. We will show that
\[
H^p(\bbA^\hbar_M
\otimes^{\bbL}_{(\bbA^\hbar_M)^{\frak{e}}}
\bbA^\hbar_M)[\hbar^{-1}]\isomo
\left\lbrace\begin{array}{lll} 0 & \text{if} & p\neq -2d \\
\bbC_M[\hbar^{-1},\hbar ]] & \text{if} & p = -2d \end{array}\right.
\]

Consider a ``local trivialization'' of $\bbA^\hbar_M$ as in 
\eqref{loc-iso-h}. It induces an isomorphism
\begin{equation}\label{loc-iso-hoch}
H^{-2d}(\bbA^\hbar_{U'}\otimes^{\bbL}_{(\bbA^\hbar_{U'})^{\frak{e}}}
\bbA^\hbar_{U'})[\hbar^{-1}] @>{\isomo}>> 
H^{-2d}(\bbA^\hbar_{U}\otimes^{\bbL}_{(\bbA^\hbar_{U})^{\frak{e}}}
\bbA^\hbar_{U})[\hbar^{-1}]\ .
\end{equation}

The expression $Alt(1\otimes x_1\otimes\cdots\otimes x_d\otimes
\displaystyle\frac{\xi_1}{\hbar}\otimes\cdots\otimes\displaystyle
\frac{\xi_d}{\hbar})$ represents a (non-trivial) global section of 
$H^{2d}(\bbA^\hbar_{U'}\otimes^{\bbL}_{(\bbA^\hbar_{U'})^{\frak{e}}}
\bbA^\hbar_{U'})[\hbar^{-1}]$ whose 
image $\Phi^{\bbA}_U$ under the isomorphism \eqref{loc-iso-hoch} is
independent of the local trivialization \eqref{loc-iso-h}.

Let
\[
\mu_{\bbA}^\hbar : \bbA^\hbar_M\otimes^{\bbL}_{(\bbA^\hbar_M)^{\frak{e}}}
\bbA^\hbar_M[\hbar^{-1}] @>>> \bbC_M[\hbar^{-1},\hbar ]][2d]
\]
denote the isomorphism in the derived category of sheaves which corresponds
to the global section of 
$H^{2d}(\bbA^\hbar_M\otimes^{\bbL}_{(\bbA^\hbar_M)^{\frak{e}}}
\bbA^\hbar_M)[\hbar^{-1}]$ determined by the condition that it
restricts to $\Phi^{\bbA}_U$ for every sufficiently small open set
$U\subset M$ as above.

We now turn to the deformation quantization $\bbA^\hbar_{T^*X}$ as in
\ref{subsec:microloc}
The compatibility of all of the maps defined above is expressed by the
commutativity of the following diagram:
\[
\begin{CD}
\pi^{-1}\cD_X @>>> \cE_X @>>> \bbA^\hbar_{T^*X}[\hbar^{-1}] \\
@VVV 		 @VVV 		 @VVV \\
\pi^{-1}\left(\cD_X\otimes^{\bbL}_{\cD_X^{\frak{e}}}\cD_X\right) @>>>
\cE_X\otimes^{\bbL}_{\cE_X^{\frak{e}}}\cE_X @>>>
\bbA^\hbar_{U}\otimes^{\bbL}_{(\bbA^\hbar_{U})^{\frak{e}}}
\bbA^\hbar_{U}[\hbar^{-1}] \\
@V{\mu^\hbar_{\cD}}VV @V{\mu^\hbar_{\cE}}VV @V{\mu^\hbar_{\bbA}}VV \\
\pi^{-1}\bbC_X[\hbar^{-1},\hbar]][2d] @>>>
\bbC_{T^*X}[\hbar^{-1},\hbar]][2d] @>>> \bbC_{T^*X}[\hbar^{-1},\hbar]][2d]
\end{CD}
\]

Although we will restrict ourselves to the discussion of the cyclic homology
of $\bbA^\hbar_M[\hbar^{-1}]$, analogs of the statements below hold for the 
algebras $\cD_X$ and $\cE_X$.

Since there are no nontirvial morphisms
$\bbC_M[\hbar^{-1},\hbar]]\to\bbC_M[\hbar^{-1},\hbar]][-1]$ in the derived
category it follows that the map
$B : C_\bullet(\bbA^\hbar_M)[\hbar^{-1}]\to
C_\bullet(\bbA^\hbar_M)[\hbar^{-1}][-1]$ represents the trivial morphism
(in the derived category) and, consequently, there are isomorphisms
(in the derived category)
\begin{equation*}
CC^-_\bullet(\bbA^\hbar_M)[\hbar^{-1}]\isomoto
\prod_{p = 0}^{\infty} C^\bullet(\bbA^\hbar_M)[\hbar^{-1}][-2p]\isomoto
\prod_{p = 0}^{\infty} \bbC_M[\hbar^{-1},\hbar]][2d - 2p]
\end{equation*}
and
\begin{equation}\label{map:TR}
CC^{per}_\bullet(\bbA^\hbar_M)[\hbar^{-1}]\isomoto
\prod_{p = -\infty}^{\infty} C^\bullet(\bbA^\hbar_M)[\hbar^{-1}][-2p]\isomoto
\prod_{p = -\infty}^{\infty} \bbC_M[\hbar^{-1},\hbar]][2d - 2p]
\end{equation}
which are induced by $\mu^\hbar_{\bbA}$ on each factor. We will denote the
latter composition by $\tilde\mu^\hbar_{\bbA}$.

It is not difficult to show that the inverse to $\tilde\mu^\hbar_{\bbA}$
is provided by the map of periodic cyclic complexes induced by the inclusion
$\bbC_M[[\hbar]]\hookrightarrow\bbA^\hbar_M$.

\subsection{Euler classes of $\cD_X$-, $\cE_X$- and $\bbA^\hbar_{T^*X}$-modules}
Consider a perfect complex $\cM^\bullet$ of $\cD_X$-modules and a closed
subvariety $\Lambda$ of $T^*X$ containing $\chr(\cM^\bullet)$. It is well
known that $\chr(\cM^\bullet) = \supp(\pi^{-1}\cM^\bullet
\otimes_{\pi^{-1}\cD_X}\cE_X)$ and that the microlocal Euler class
$\mu\eu(\cM^\bullet)\in H^{2d}_\Lambda(T^*X;\bbC)$ depends only
on the microlocalization $\pi^{-1}\cM^\bullet\otimes_{\pi^{-1}\cD_X}\cE_X$.
In fact, it is not difficult to establish the equality
\[
\mu\eu(\cM^\bullet) = \mu_{\cE}\left(\eu^\Lambda_{\cE}
(\pi^{-1}\cM^\bullet\otimes_{\pi^{-1}\cD_X}\cE_X)\right)\ .
\]

Consider a perfect complex $\cN^\bullet$ of $\cE_X$-modules and a closed
subset $\Lambda$ of $T^*X$ containing $\supp(\cN^\bullet)$. The commutativity
of the diagram
\[
\begin{CD}
\R\shHom_{\cE_X}(\cN^\bullet,\cN^\bullet) & @>>> &
\R\shHom_{\bbA^\hbar_{T^*X}[\hbar^{-1}]}(\cN^\bullet\otimes_{\cE_X}
\bbA^\hbar_{T^*X}[\hbar^{-1}],\cN^\bullet\otimes_{\cE_X}
\bbA^\hbar_{T^*X}[\hbar^{-1}]) \\
@V{\cL_{\cE}}VV & & 
@VV{\cL_{\bbA^\hbar_{T^*X}[\hbar^{-1}]}}V \\
\cE\otimes^{\bbL}_{\cE_X^{\frak{e}}}\cE_X &
@>>> & \bbA^\hbar_{T^*X}\otimes^{\bbL}_{(\bbA^\hbar_{T^*X})^{\frak{e}}}
\bbA^\hbar_{T^*X}[\hbar^{-1}]
\end{CD}
\]
implies the identity $\mu^\hbar_{\cE}\circ\Eu_{\cE}(\cN^\bullet) =
\mu^\hbar_{\bbA}\circ\Eu_{\bbA^\hbar_{T^*X}[\hbar^{-1}]}
(\cN^\bullet\otimes_{\cE_X}\bbA^\hbar_{T^*X}[\hbar^{-1}])$, and,
consequently, the equality
\[
\mu^\hbar_{\cE}\left(\eu^\Lambda_{\cE}(\cN^\bullet)\right) =
\mu^\hbar_{\bbA}\left(\eu^\Lambda_{\bbA}
(\cN^\bullet\otimes_{\cE_X}\bbA^\hbar_{T^*X}[\hbar^{-1}])\right)\ .
\]

Thus, calculation of microlocal Euler classes reduces to calculation
of Euler classes for $\bbA^\hbar_M[\hbar^{-1}]$-modules.

\section{Riemann-Roch type theorems}
\label{section:RR}
\subsection{The Riemann-Roch theorem for periodic cyclic cocycles}
The following theorem constitutes the central result of this note.

\begin{thm}\label{thm:main}
The diagram (in the derived category of Abelian sheaves on $M$)
\[
\begin{CD}
CC^{per}_\bullet(\bbA^\hbar_M) @>{\sigma}>> CC^{per}_\bullet(\cO_M) \\
@V{\iota}VV @VV
{\tilde\mu_{\cO}\smile \widehat A(TM)\smile e^\theta}V \\
CC^{per}_\bullet(\bbA^\hbar_M)[\hbar^{-1}] @>{\tilde\mu^\hbar_{\bbA}}>>
\displaystyle\prod_{p=-\infty}^\infty \bbC_M[\hbar^{-1},\hbar]][-2p]
\end{CD}
\]
is commutative.
\end{thm}

The proof of Theorem \ref{thm:main} is postponed until the later sections.
In Section \ref{section:GF} we introduce the methods of Gel'fand-Fuchs
cohomology and reduce (see Corollary \ref{cor:reduction})
Theorem \ref{thm:main} to an analogous statement (Theorem \ref{thm:formalRR})
in the case when $M$ is a formal neighborhood of the origin in a symplectic
vector space over $\bbC$ which is formulated in Section \ref{section:formalRR}.
The rest of this section is devoted to corollaries of Theorem \ref{thm:main}.

\subsection{Riemann-Roch for $\bbA^\hbar$-modules}
The commutativity of the diagram
\[
\begin{CD}
K^0_\Lambda(\bbA^\hbar_M[\hbar^{-1}]) @<{\iota}<< K^0_\Lambda(\bbA^\hbar_M) 
@>{\sigma}>> K^0_\Lambda(\cO_M) \\
@V{ch^\Lambda_{\bbA^\hbar[\hbar^{-1}]}}VV @V{ch^\Lambda_{\bbA^\hbar}}VV
@VV{ch^\Lambda_{\cO}}V \\
H^0_\Lambda(M;CC^{per}_\bullet(\bbA^\hbar_M)[\hbar^{-1}]) @<{\iota}<<
H^0_\Lambda(M;CC^{per}_\bullet(\bbA^\hbar_M))
@>{\sigma}>> H^0_\Lambda(M;CC^{per}_\bullet(\cO_M))
\end{CD}
\]
yields the following.

\begin{cor}\label{cor:rrAch}
Suppose that $M$ is a complex manifold, $\bbA^\hbar_M$ is a symplectic
deformation quantization of $M$, $\cM^\bullet$ is a perfect complex of
$\bbA^\hbar_M$-modules and $\Lambda$
is a closed subvariety of $M$ containing $\supp(\cM^\bullet)$. Then
\begin{equation}\label{eq:rrAch}
\mu^\hbar_{\bbA}\left(ch^\Lambda_{\bbA^\hbar_M[\hbar^{-1}]}
(\cM^\bullet[\hbar^{-1}])\right)
= \tilde\mu_{\cO}(ch^\Lambda_{\cO}(\sigma(\cM^\bullet)))
\smile\widehat A(TM)\smile e^\theta
\end{equation}
in $H_\Lambda^\bullet(M;\bbC[\hbar^{-1},\hbar ]])$,
where $\theta$ is the characteristic class
of the deformation quantization $\bbA^\hbar_M$.
\end{cor}

Note that the class $\widehat A(E)$ is defined for any symplectic vector
bundle $E$, for example by choosing a reduction of the (symplectic) structure
group of $E$ to the unitary group. 

Recall that, for an element $\alpha$ of a graded object, $[\alpha]^p$ denotes
the homogeneous component of $\alpha$ of degree $p$.

\begin{cor}
Under the assumptions of Corollary \ref{cor:rrAch}
\begin{equation}\label{eq:rrAeu}
\mu^\hbar_{\bbA}\left(\eu^\Lambda_{\bbA^\hbar_M[\hbar^{-1}]}
(\cM^\bullet[\hbar^{-1}])\right)
= \left[\tilde\mu_{\cO}(ch^\Lambda_{\cO}(\sigma(\cM^\bullet)))
\smile\widehat A(TM)\smile e^\theta)
\right]^{\dim_{\bbC}M}
\end{equation}
in $H_\Lambda^\bullet(M;\bbC[\hbar^{-1},\hbar ]])$.
\end{cor}

\subsection{Riemann-Roch for $\cD$- and $\cE$-modules}
If $M = T^*X$ for a complex manifold $X$, and $\bbA^\hbar_{T^*X}$
is the deformation quantization with the characteristic class
$\theta = \frac12\pi^*c_1(X)$, then
\[
\widehat A(TM)\smile e^\theta = \pi^*Td(TX)
\]
and the right hand side of \eqref{eq:rrAeu} is, clearly, independent of $\hbar$;
if $\cM^\bullet[\hbar^{-1}]$ is obtained by an extension of scalars from
a complex of $\cE_X$-modules, then, clearly, so is the left hand side.

Thus, we obtain Conjecture \ref{conj:ch} of P.Schapira and J.-P.Schneiders.
\begin{cor}
Suppose that $X$ is a complex manifold, $(\cM^\bullet,F_\bullet)$ is a perfect 
complex of $\cD_X$-modules with a good filtration and $\Lambda$ is a closed
subvariety of $T^*X$ containing $\chr\cM^\bullet$. Then
\[
\mu\eu_\Lambda(\cM^\bullet) =
\left[ ch^\Lambda_{\cO_{T^*X}}(\sigma(\cM^\bullet))
\smile\pi^*Td(TX)\right]^{2\dim_{\bbC}X}
\]
in $H_\Lambda^{2\dim_{\bbC}X}(T^*X;\bbC)$.
\end{cor}

\section{The Riemann-Roch formula in the formal setting}
\label{section:formalRR}
The Weyl algebra of a symplectic vector space $(V,\omega)$ over $\bbC$
may be considered as a symplectic deformation quantization of
the completion of $V$ at the origin. In this section we introduce the
notations and the facts necessary to state the analogue of Theorem
\ref{thm:main} in this setting.

In what follows $(V,\omega)$ is viewed as a symplectic manifold. 

\subsection{The Weyl algebra}
Here we briefly recall the definition and the basic properties of the
Weyl algebra $W = W(V)$ of a (finite dimensional) symplectic vector space 
$(V,\omega)$ over $\bbC$. Let $V^\ast = \Hom_{\bbC}(V,\bbC)$.

Let $I=I(V)$ denote the kernel of the map (of $\bbC$-algebras)
\[
\Symm^\bullet (V^\ast)\otimes\bbC[\hbar] @>>> \bbC\ .
\]
Let
\[
\widehat\Symm^\bullet(V^\ast)[[\hbar]] = \varprojlim
\frac{\Symm^\bullet (V^\ast)\otimes\bbC[\hbar]}{I^n}
\]
and let $\widehat I = \widehat I(V)$ denote the kernel of the map
$\widehat\Symm^\bullet(V^\ast)[[\hbar]]\to\bbC$.

The Moyal-Weyl product on $\widehat\Symm^\bullet(V^\ast)[[\hbar]]$ is defined
by the formula
\begin{equation}\label{formula:MW}
f\ast g = \sum_{n=0}^\infty\frac{1}{n!}\left(\frac{\sqrt{-1}\hbar}{2}\right)^n
\omega(d^nf,d^ng)\ ,
\end{equation}
where
\[
d^n : \widehat\Symm^\bullet(V^\ast) @>>>
\widehat\Symm^\bullet(V^\ast)\otimes \Symm^n(V)
\]
assigns to a jet of a function the symmetric tensor composed of its
$n$-th order partial deriviatives, and $\omega$ is extended naturally to
a bilinear form on $\Symm^n(V)$.

The Moyal-Weyl product
endows $\widehat\Symm^\bullet(V^\ast)[[\hbar]]$ with a structure of an
associative algebra with unit over $\bbC[[\hbar]]$ which contains
$\widehat I$ as a twosided ideal. Moreover, the Moyal-Weyl product is
continuous in the $\widehat I$-adic topology. Let $W=W(V)$ denote
the topological algebra over $\bbC[[\hbar]]$ whose underlying
$\bbC[[\hbar]]$-module is $\widehat\Symm^\bullet(V^\ast)[[\hbar]]$ and the
multiplication is given by the Moyal-Weyl product.

Let $F_pW = \widehat I^{-p}$. Then $(W,F_\bullet)$ is a filtered ring.
Note also that the center of $W$ is equal to $\bbC[[\hbar]]$ and
$\left[W,W\right] = \hbar\cdot W$.

Clearly, the association $(V,\omega)\mapsto W(V)$ is functorial. In
particular, the group $Sp(V)$ acts naturally on $W(V)$ by continuous
algebra automorphisms.

\subsection{Derivations of the Weyl algebra}
Let $\frak{g} = \frak{g}(V)$ denote the Lie algebra of continuous,
$\bbC[[\hbar]]$-linear derivations of $W$. Then there is a central extension
of Lie algebras
\[
0 @>>> \frac{1}{\sqrt{-1}\hbar}\bbC[[\hbar]] @>>>
\frac{1}{\sqrt{-1}\hbar}W @>>> \frak{g} @>>> 0\ ,
\]
where the Lie algebra structure on $\frac{1}{\sqrt{-1}\hbar}W$ is
given by the commutator (note that
$\left[\frac{1}{\sqrt{-1}\hbar}W,\frac{1}{\sqrt{-1}\hbar}W\right]
\subseteq\frac{1}{\sqrt{-1}\hbar}W$) and the second map is defined
by $\frac{1}{\sqrt{-1}\hbar}f\mapsto \frac{1}{\sqrt{-1}\hbar}\left[f,\bullet
\right]$.

Let
\[
F_p\frak{g} = \left\{D\in\frak{g}\ \vert\ D(F_iW)\subseteq F_{i+p}W
\ \text{for all $i$}\right\}
\]
Then $(\frak{g},F_\bullet)$ is a filtered Lie algebra and the action
of $\frak{g}$ on $W$ respects the filtrations, i.e.
$\left[F_p\frak{g},F_q\frak{g}\right]\subseteq F_{p+q}\frak{g}$ and
$F_p\frak{g}F_qW\subseteq F_{p+q}W$.

The following properties of the filtered Lie algebra $(\frak{g},F_\bullet)$
are easily verified:
\begin{enumerate}
\item $Gr^F_p\frak{g} = 0$ for $p>1$ (in particular $\frak{g}=F_1\frak{g}$),
hence $Gr^F_1\frak{g}$ is Abelian;

\item the composition $\frac{1}{\sqrt{-1}\hbar}V\hookrightarrow
\frac{1}{\sqrt{-1}\hbar}W\to\frak{g}\to Gr^F_1\frak{g}$ is an
isomorphism;

\item the composition $\frak{sp}(V)\to\Symm^2(V^\ast)\to
\frac{1}{\sqrt{-1}\hbar}\widehat I\to F_0\frak{g}\to Gr^F_0\frak{g}$
is an isomophism;

\item under the above isomorphisms the action of $Gr^F_0\frak{g}$ on
$Gr^F_1\frak{g}$ is identified with the natural action of
$\frak{sp}(V)$ on $V$ (in particular there is an isomorphism
$\frak{g}/F_{-1}\frak{g}\isomo V\ltimes\frak{sp}(V)$);

\item the Lie algebra $F_{-1}\frak{g}$ is pro-nilpotent.
\end{enumerate}

In what follows $\frak{h}$ will denote the Lie subalgebra of $\frak{g}$
which is the image of the embedding $\frak{sp}(V)\hookrightarrow\frak{g}$.

\subsection{The Weyl algebra as a deformation quantization}
The (commutative) algebra $W/\hbar\cdot W$ is naturally isomorphic to
the completion $\widehat\cO =\widehat\cO_V$ of the ring $\cO_V$ of regular functions
on $V$ at the origin, i.e. with respect to the powers of the maximal ideal
$\frak{m}$ of functions which vanish at $0\in V$. The natural surjective
map
\[
\sigma : W(V) @>>> \widehat\cO_V
\]
is strictly compatible with the $\widehat I$-adic filtration on $W$ and the
$\frak{m}$-adic filtration on $\widehat\cO_V$.

The Lie algebra $\frak{g}$ acts by derivations on $\widehat\cO$
by the formula
\[
D(f) = \sigma(D(\tilde f))\ ,
\]
where $D\in \frak{g}$ and $\tilde f\in W$ is such that $\sigma(\tilde f)=f$.
Thus, $\sigma$ is a map of $\frak{g}$-modules.

\subsection{The Hochschild homology of the Weyl algebra}
We recall the calculation of the Hochschild homology of the Weyl algebra
(\cite{FT1}, \cite{Bry}).

The Hochschild homology of $W$ may be computed using the Koszul
resolution of $W$ as a $W^{\frak{e}}\overset{def}{=}
W\widehat\otimes_{\bbC[[\hbar]]}W^{op}$ module.

The Koszul complex
$(K^\bullet,\partial)$ is defined by
\[
K^{-q} = W\otimes{\bigwedge}^q V^\ast\otimes W
\]
with the differential acting by
\begin{eqnarray*}
\partial(f\otimes v_1\wedge\ldots\wedge v_q\otimes g) & = &
\sum_i (-1)^i fv_i\otimes v_1\wedge\ldots\wedge\widehat v_i\wedge\ldots\wedge
v_q \\
  & + & \sum_i (-1)^i f\otimes v_1\wedge\ldots\wedge\widehat v_i\wedge\ldots
  \wedge v_q\otimes v_ig\ .
\end{eqnarray*}
Here we consider $V^\ast$ embedded in $W$ and $\bigwedge^qV^\ast$ embedded
in $(V^\ast)^{\otimes q}$.

The map $K^\bullet\to W$ of complexes of $W^{\frak{e}}$-modules with the
only nontrivial component (in degree zero) given by multiplication is easily
seen to be a quasiisomorphism. The map
\[
K^\bullet(W)\overset{def}{=}K^\bullet\otimes_{W^{\frak{e}}}W @>>> C_\bullet(W)
\]
defined by
\[
f\otimes v_1\wedge\ldots\wedge\ldots\wedge v_q\otimes g\otimes h
\mapsto fhg\otimes Alt(v_1\otimes\cdots\otimes v_q)
\]
is easily seen to be a quasiisomorphism. Hence it induces a quasiisomorphism
\[
K^\bullet(W)\otimes_{\bbC[[\hbar]]}\bbC[\hbar^{-1},\hbar]]
@>>> C_\bullet(W)\otimes_{\bbC[[\hbar]]}\bbC[\hbar^{-1},\hbar]]\ .
\]
Below we will use $(\bullet)[\hbar^{-1}]$ to denote
$(\bullet)\otimes_{\bbC[[\hbar]]}\bbC[\hbar^{-1},\hbar]]$.

Under the natural isomorphism
\[
K^{-q}(W) @>>> W\otimes{\bigwedge}^qV^\ast
\]
defined by
\[
f\otimes v_1\wedge\ldots\wedge v_q\otimes g\otimes h\mapsto fhg\otimes
v_1\wedge\ldots\wedge v_q
\]
the induced differential acts on the latter by
\[
\partial(f\otimes v_1\wedge\ldots\wedge v_q) =
\sum_i [f,v_i]\otimes v_1\wedge\ldots\wedge\widehat v_i\wedge\ldots\wedge v_q\ .
\] 

Suppose that $\dim V = 2d$.
Let $\widehat\Omega^\bullet =\widehat\Omega^\bullet_V$ denote the de Rham
complex of $V$ with
formal coefficients (i.e. $\widehat\Omega^q_V =\Omega^q_V\otimes_{\cO_V}
\widehat\cO_V$). 
The map
\[
W\otimes{\bigwedge}^q V^\ast @>>>
\widehat\Omega^{2d-q}[[\hbar]]
\]
given, by
\[
f\otimes v_1\wedge\ldots\wedge v_q \mapsto
f\cdot\iota_{v_1}\cdots\iota_{v_q}(\omega^{\wedge d})
\]
(where $f\in\widehat\cO$)
is easily seen to determine an isomorphism of complexes
\[
(K^\bullet(W),\partial)
@>>> (\widehat\Omega^\bullet[[\hbar]], \hbar\cdot d)[2d]\ .
\]

The formal Poincare Lemma implies that
\[
H^q(K^\bullet\otimes_{W^{\frak{e}}}W[\hbar^{-1}]) \isomo
HH_q(W;W)[\hbar^{-1}]\isomo\left\lbrace\begin{array}{ll}
0 & \text{if $q\neq 2d$} \\ \bbC[\hbar^{-1},\hbar]] & \text{if $q=2d$}
\end{array}\right.
\]
The canonical generator of $HH_{2d}(W;W)[\hbar^{-1}]$ is represented by
the cycle $1\otimes\frac{1}{d!}\omega^{\wedge d}\in W\otimes
{\bigwedge}^{2d}V^\ast$.
If $x_1,\ldots ,x_d,\xi_1,\ldots ,\xi_d$ is (dual to) a symplectic basis
of $V$ (so that $\omega = \sum_i x_i\wedge\xi_i$), then the canonical generator
is represented by the cycle $Alt(1\otimes x_1\otimes\cdots\otimes x_d\otimes
\xi_1\otimes\cdots\otimes\xi_d)\in C_{2d}(W)$. We will denote this cycle
(and its class) by $\Phi = \Phi_V$. Note also that $\Phi$ corresponds
under the above (quasi)isomorphisms to the cocycle $1\in\widehat\Omega^0$.

Observe that the Lie algebra $\frak{g}$ acts on all of the complexes
introduced above (by Lie deriviative) and all maps defined above are, in fact,
$\frak{g}$-equivariant. The cycle $\Phi$ is not invariant under the action
of $\frak{g}$. It is, however, invariant under the action of the subalgebra
$\frak{h}$.

\subsection{Characteristic classes in Lie algebra cohomology}
We will presently construct the classes in relative Lie algebra cohomology
of the pair $(\frak{g},\frak{h})$ which enter the Riemann-Roch formula in
the present setting.

\subsubsection{The trace density}
Since $\frak{h}$ acts semi-simply on $C_\bullet(W)[\hbar^{-1}]$ and
$\widehat\Omega^\bullet[\hbar^{-1},\hbar]]$, the quasiisomorphism
\[
\widehat\Omega^\bullet[\hbar^{-1},\hbar]][2d] @>>> C_\bullet(W)[\hbar^{-1}]
\]
constructed above admits an $\frak{h}$-equivariant splitting
\[
\mu^\hbar_{(0)} :
C_\bullet(W)[\hbar^{-1}] @>>> \widehat\Omega^\bullet_V[\hbar^{-1},\hbar]][2d]
\]
which is a quasiisomorphism.

We will consider the map $\mu^\hbar_{(0)}$ as a relative Lie algebra cochain
\[
\mu^\hbar_{(0)}\in
C^0(\frak{g},\frak{h};\Hom^0(C_\bullet(W)[\hbar^{-1}],
\widehat\Omega^\bullet[\hbar^{-1},\hbar]][2d]))\ .
\]

\begin{lemma}
$\mu^\hbar_{(0)}$ extends to a cocycle $\mu^\hbar = \sum_p\mu^\hbar_{(p)}$
with
\[
\mu^\hbar_{(p)}\in
C^p(\frak{g},\frak{h};\Hom^{-p}(C_\bullet(W)[\hbar^{-1}],
\widehat\Omega^\bullet[\hbar^{-1},\hbar]][2d]))\ .
\]
Moreover, any two such extensions are cohomologous.
\end{lemma}

Since the complex
$C^\bullet(\frak{g},\frak{h};\Hom^\bullet_{\bbC[[\hbar]]}
(C_\bullet(W)[\hbar^{-1}],\widehat\Omega^\bullet[\hbar^{-1},\hbar]][2d]))$
represents (in the derived category) the object
$\R\Hom^\bullet_{(\frak{g},\frak{h})}(C_\bullet(W)[\hbar^{-1}],
\widehat\Omega^\bullet[\hbar^{-1},\hbar]][2d])$, $\mu^\hbar$ represents
a well defined isomorphism
\[
\mu^\hbar :
C_\bullet(W)[\hbar^{-1}] @>>> \widehat\Omega^\bullet[\hbar^{-1},\hbar]][2d]
\]
in the derived category of $(\frak{g},\frak{h})$-modules. The image of
$\mu^\hbar$ under the functor of forgetting the module structure is
$\mu^\hbar_{(0)}$.

Cup product with $\mu^\hbar$ induces the quasiisomorphism of complexes
\[
\mu^\hbar : C^\bullet(\frak{g},\frak{h};C_\bullet(W)[\hbar^{-1}]) @>>>
C^\bullet(\frak{g},\frak{h};\widehat\Omega^\bullet[\hbar^{-1},\hbar]][2d])
\]
unique up to homotopy.

\begin{lemma}
$\mu^\hbar_{(0)}$ extends to an $\frak{h}$-equivariant quasiisomorphism of
complexes
\[
\tilde\mu^\hbar_{(0)} : CC^{per}_\bullet(W)[\hbar^{-1}] @>>>
\prod_{p\in\bbZ}\widehat\Omega^\bullet[\hbar^{-1},\hbar]][2p]\ .
\]
\end{lemma}

We will consider the map $\tilde\mu^\hbar_{(0)}$ as a relative Lie algebra
cochain
\[
\tilde\mu^\hbar_{(0)}\in
C^0(\frak{g},\frak{h};\Hom^0(CC^{per}_\bullet(W)[\hbar^{-1}],
\prod_{p\in\bbZ}\widehat\Omega^\bullet[\hbar^{-1},\hbar]][2p]))\ .
\]

\begin{lemma}\label{lemma:formalTR}
$\tilde\mu^\hbar_{(0)}$ extends to a cocycle $\tilde\mu^\hbar =
\sum_p\tilde\mu^\hbar_{(p)}$
with
\[
\mu^\hbar_{(p)}\in
C^p(\frak{g},\frak{h};\Hom^{-p}(CC^{per}_\bullet(W)[\hbar^{-1}],
\prod_{q\in\bbZ}\widehat\Omega^\bullet[\hbar^{-1},\hbar]][2q]))\ .
\]
Moreover, any two such extensions are cohomologous.
\end{lemma}

Since the complex
$C^\bullet(\frak{g},\frak{h};\Hom^\bullet
(CC^{per}_\bullet(W)[\hbar^{-1}],\prod_{p\in\bbZ}
\widehat\Omega^\bullet[\hbar^{-1},\hbar]][2p]))$ represents the object
$\R\Hom^\bullet_{(\frak{g},\frak{h})}(CC^{per}_\bullet(W)[\hbar^{-1}],
\prod_{p\in\bbZ}\widehat\Omega^\bullet[\hbar^{-1},\hbar]][2d])$,
$\tilde\mu^\hbar$ represents a well defined isomorphism
\[
\tilde\mu^\hbar :
CC^{per}_\bullet(W)[\hbar^{-1}] @>>>
\prod_{p\in\bbZ}\widehat\Omega^\bullet[\hbar^{-1},\hbar]][2d]
\]
in the derived category of $(\frak{g},\frak{h})$-modules. The image of
$\tilde\mu^\hbar$ under the functor of forgetting the module structure
is $\tilde\mu^\hbar_{(0)}$.

Cup product with $\tilde\mu^\hbar$ induces the quasiisomorphism of complexes
\[
\tilde\mu^\hbar : C^\bullet(\frak{g},\frak{h};CC^{per}_\bullet(W)[\hbar^{-1}])
@>>>
C^\bullet(\frak{g},\frak{h};\prod_{p\in\bbZ}
\widehat\Omega^\bullet[\hbar^{-1},\hbar]][2p])
\]
unique up to homotopy.

The natural inclusion
\[
\iota : CC^{per}_\bullet(W) @>>> CC^{per}_\bullet(W)[\hbar^{-1}]
\]
is a morphism of complexes of $(\frak{g},\frak{h})$-modules, therefore
determines a cocycle
\[
\iota\in C^0(\frak{g},\frak{h};\Hom^0
(CC^{per}_\bullet(W), CC^{per}_\bullet(W)[\hbar^{-1}]))\ .
\]

The cup product of $\iota$ and $\tilde\mu^\hbar$ is a cocycle
\[
\tilde\mu^\hbar\smile\iota\in
C^\bullet(\frak{g},\frak{h};\Hom^\bullet
(CC^{per}_\bullet(W),\prod_{p\in\bbZ}
\widehat\Omega^\bullet[\hbar^{-1},\hbar]][2p])
\]
of (total) degree zero which prepresents the morphism
\[
\tilde\mu^\hbar\circ\iota :
CC^{per}_\bullet(W) @>>> \prod_{p\in\bbZ}
\widehat\Omega^\bullet[\hbar^{-1},\hbar]][2p]
\]
in the derived category of $(\frak{g},\frak{h})$-modules.

\subsubsection{The symbol and the Hochschild-Kostant-Rosenberg map}
\label{sssection:HKR}
The maps
\[
\sigma : CC^{per}_\bullet(W) @>>> CC^{per}_\bullet(\widehat\cO)
\]
(induced by $\sigma : W\to\widehat\cO$) and
\[
\tilde\mu : CC^{per}_\bullet(\widehat\cO) @>>> \prod_{p\in\bbZ}
\widehat\Omega^\bullet[2p] @>>> \prod_{p\in\bbZ}
\widehat\Omega^\bullet[2p][\hbar^{-1},\hbar]]
\]
(defined by $f_0\otimes\cdots\otimes f_p\mapsto\frac{1}{p!}
f_0df_1\wedge\ldots\wedge df_p$)
are morphisms of complexes of $(\frak{g},\frak{h})$-modules, therefore
determine cocycles
\[
\sigma\in C^\bullet(\frak{g},\frak{h};\Hom^\bullet
(CC^{per}_\bullet(W), CC^{per}_\bullet(\widehat\cO)))
\]
and
\[
\tilde\mu\in C^\bullet(\frak{g},\frak{h};\Hom^\bullet
(CC^{per}_\bullet(\widehat\cO),\prod_{p\in\bbZ}
\widehat\Omega^\bullet [\hbar^{-1},\hbar]][2p]))\ .
\]

\subsubsection{The characteristic class of the deformation}
\label{sssection:charcl}
The central extension of Lie algebras
\[
0 @>>> \frac{1}{\sqrt{-1}\hbar}\bbC[[\hbar]] @>>>
\frac{1}{\sqrt{-1}\hbar}W @>>> \frak{g} @>>> 0
\]
restricts to a trivial extension of $\frak{h}$, therefore is classified by
a class $\theta\in H^2(\frak{g},\frak{h};\frac{1}{\sqrt{-1}\hbar}\bbC[[\hbar]])$
represented by the cocycle
\[
\theta : X\wedge Y\mapsto \widetilde{[X,Y]} -
[\widetilde{X},\widetilde{Y}]
\]
where $\widetilde{(\ )}$ is a choice of a $\bbC[[\hbar]]$-linear splitting
of the extension.

\subsubsection{The $\widehat A$-class}\label{sssection:formalA}
Let $\nabla : \frak{g}\to\frak{h}$ denote an $\frak{h}$-equivariant splitting
of the inclusion $\frak{h}\hookrightarrow\frak{g}$ and let
\[
R(X,Y) = [\nabla (X),\nabla (Y)]-\nabla ([X,Y])
\]
for $X,\ Y\in\frak{g}$. Then, $R\in C^2(\frak{g},\frak{h};\frak{h})$ is
a cocycle (where $\frak{h}$ is considered as a trivial
$(\frak{g},\frak{h})$-module). The $k$-fold cup product of $R$ with itself
is a cocycle $R^{\smile k}\in C^{2k}(\frak{g},\frak{h};\frak{h}^{\otimes k})$.

The splitting $\nabla$ determines the Chern-Weil map
\[
CW : \widehat\Symm^\bullet(\frak{h}^\ast)^{\frak{h}} @>>>
C^{2\bullet}(\frak{g},\frak{h};\bbC)
\]
of complexes (with $\widehat\Symm^\bullet(\frak{h}^\ast)^{\frak{h}}$ endowed
with the trivial differential) which is defined as the composite
\[
\Symm^k(\frak{h}^\ast)^{\frak{h}} @>>>
C^0(\frak{h};\Symm^k(\frak{h}^\ast)) @>>>
C^0(\frak{h};\Symm^k(\frak{h}^\ast))\otimes C^{2k}(\frak{g},\frak{h};
\frak{h}^{\otimes k}) @>>>
C^{2k}(\frak{g},\frak{h};\bbC)
\]
where the first map is the natural inclusion (of the cocycles), the second map
is $P\mapsto P\otimes R^{\smile k}$, and the last map is induced by the
natural pairing on the coefficients.

Let $\widehat A$ denote the image under the Chern-Weil map of
\[
\frak{h}\ni X\mapsto\det\left(\frac{ad(\frac{X}{2})}{exp(ad(\frac{X}{2})) -
exp(ad(-\frac{X}{2}))}\right)\ .
\]

\subsection{Riemann-Roch formula in Lie algebra cohomology}
The following theorem is the analog of the Theorem \ref{thm:main} in the
setting of this section.

\begin{thm}\label{thm:formalRR}
The cocycle $\tilde\mu^\hbar\smile\iota - \tilde\mu\smile\widehat A\smile
e^\theta\smile\sigma$ is cohomolgous to zero in
$C^\bullet(\frak{g},\frak{h};\Hom^\bullet
(CC^{per}_\bullet(W),\prod_{p\in\bbZ}
\widehat\Omega^\bullet [\hbar^{-1},\hbar]][2p]))$.
\end{thm}

\begin{cor}
The diagram
\[
\begin{CD}
CC^{per}_\bullet(W) @>{\sigma}>> CC^{per}_\bullet(\widehat\cO) \\
@V{\iota}VV @VV{\tilde\mu\smile\widehat A\smile e^\theta}V \\
CC^{per}_\bullet(W)[\hbar^{-1}] @>{\tilde\mu^\hbar}>>
\prod_{p\in\bbZ}\widehat\Omega^\bullet_V[\hbar^{-1},\hbar]][2p]
\end{CD}
\]
in the derived category of $(\frak{g},\frak{h})$-modules is commutative.
\end{cor}

\begin{cor}
The diagram
\[
\begin{CD}
C^\bullet(\frak{g},\frak{h};CC^{per}_\bullet(W)) @>{\sigma}>>
C^\bullet(\frak{g},\frak{h};CC^{per}_\bullet(\widehat\cO)) \\
@V{\iota}VV
@VV{\tilde\mu\smile\widehat A\smile e^\theta}V
\\
C^\bullet(\frak{g},\frak{h};CC^{per}_\bullet(W)[\hbar^{-1}])
@>{\tilde\mu^\hbar}>> C^\bullet(\frak{g},\frak{h};
\prod_{p\in\bbZ}\widehat\Omega^\bullet [\hbar^{-1},\hbar]][2p])
\end{CD}
\]
is homotopy commutative.
\end{cor}

A proof of Theorem \ref{thm:formalRR} may be found in \cite{NT1}, \cite{NT2}.
In Section \ref{section:GF} we will show
that Theorem \ref{thm:main} reduces to Theorem \ref{thm:formalRR}.

\section{Gel'fand-Fuchs cohomology}\label{section:GF}
In this section we introduce the machinery of Fedosov connections and
Gel'fand-Fuchs cohomology and reduce Theorem \ref{thm:main} to the analogous
statement in the particular case when $M$ is the formal neighborhood of the
origin in a symplectic vector space over $\bbC$.

Suppose given a complex manifold $M$ and a symplectic deformation
quantization $\bbA^\hbar_M$ of $M$. Let $\omega\in H^0(M;\Omega^2_M)$
denote the associated symplectic form.

\subsection{The sheaf of Weyl algebras}
The sheaf of Weyl algebras $\bbW_M$ on $M$ is the sheaf of topological
algebras over the sheaf of topological algebras $\cO_M[[\hbar]]$ (equipped
with the $\hbar$-adic topology) defined as follows.

Let $\Theta_M$ denote the sheaf of holomorphic vector fields on $M$.
Let $\cI$ denote the kernel of the augmentation map
\[
\Symm_{\cO_M}^\bullet(\Theta_M)\otimes_{\cO_M}\cO_M[\hbar]
@>>> \cO_M\ .
\]
The completion $\widehat\Symm_{\cO_M}^\bullet(\Theta_M)[[\hbar]]$
of $\Symm_{\cO_M}^\bullet(\Theta_M)\otimes_{\cO_M}\cO_M[\hbar]$
in the $\cI$-adic topology is a topological $\cO_M[[\hbar]]$-module. The
Weyl multiplication on $\widehat\Symm_{\cO_M}^\bullet(\Theta_M)[[\hbar]]$
is defined by \eqref{formula:MW}.
Then $\bbW_M$ is the sheaf of $\cO_M[[\hbar]]$-algebras whose underlying
sheaf of $\cO_M[[\hbar]]$-modules is 
$\widehat\Symm_{\cO_M}^\bullet(\Theta_M)[[\hbar]]$ and the multiplication
is given by the Moyal-Weyl product. Note that the center of $\bbW_M$ is 
$\cO_M[[\hbar]]$, and $\left[\bbW_M,\bbW_M\right] = \hbar\cdot\bbW_M$.

Let $\widehat\cI$ denote the kernel of the canonical map
$\bbW_M\to\cO_M$. The Weyl multiplication is continuous in the
$\widehat\cI$-adic topology. 

Let $F_p\bbW_M = \widehat\cI^{-p}$. Then $(\bbW_M,F_\bullet)$ is
a filtered ring, i.e. $F_p\bbW_M\cdot F_q\bbW_M\subseteq F_{p+q}\bbW_M$.
Note that the quotients $F_p\bbW_M/F_q\bbW_M$ are locally free
$\cO_M$-modules of finite rank.

\subsection{Review of the Fedosov construction}
We refer the reader to \cite{F} and \cite{NT3} for a detailed exposition of
the construction of deformation quantizations via Fedosov connections.

Let $\cA^{p,q}_M$ denote the sheaf of complex valued $C^\infty$-forms of
type $(p,q)$ on $M$, $\cA^r_M = \oplus_{p+q=s}\cA^{p,q}_M$. Let
$F_\bullet\cA^\bullet_M$ denote the Hodge filtration and $d$ the de Rham
differential. Then $((\cA^\bullet_M,d),F_\bullet)$ is a filtered differential
graded algebra (i.e. $F_r\cA^\bullet_MF_s\cA^\bullet_M\subseteq
F_{r+s}\cA^\bullet_M$ and $d(F_r\cA^\bullet_M)\subseteq F_r\cA^\bullet_M$).

Let $\cA^p_M(\bbW_M)=\cA^p_M\otimes_{\cO_M}\bbW_M$.
Then
\[
\cA^\bullet_M(\bbW_M)\overset{def}{=}\bigoplus_p \cA^p_M(\bbW_M)[-p]
\]
has a natural structure of a sheaf of graded algebras. Let
\[
F_r\cA^\bullet_M(\bbW_M) = \sum_{p+q=r}F_p\cA^\bullet_M
\otimes_{\cO_M}F_p\bbW_M\ .
\]
Then $(\cA^\bullet_M(\bbW_M),F_\bullet)$ is a filtered graded algebra over
$(\cA^\bullet_M,F_\bullet)$ (i.e.\linebreak
$F_r\cA^p_M(\bbW_M)F_s\cA^q_M(\bbW_M)\subseteq
F_{r+s}\cA^{p+q}_M(\bbW_M)$ and $F_r\cA^p_MF_s\cA^q_M(\bbW_M)\subseteq
F_{r+s}\cA^{p+q}_M(\bbW_M)$).

Let $F_p\bbA^\hbar_M = \hbar^{-p}\cdot\bbA^\hbar_M$.

One can show that there exists a map
\[
\nabla : \cA^\bullet_M(\bbW_M) @>>> \cA^\bullet_M(\bbW_M)[1]
\]
which has the following properties:
\begin{enumerate}
\item $\nabla\left(F_p\cA^\bullet_M(\bbW_M)\right)\subseteq
F_p\cA^\bullet_M(\bbW_M)$; the induced maps
\[
Gr^F_p\nabla : Gr^F_p\cA^\bullet_M(\bbW_M) @>>> 
Gr^F_p\cA^\bullet_M(\bbW_M)[1]
\]
are $\cO_M[[\hbar]]$-linear differential operators of order one;

\item $\nabla^2 = 0$;

\item $(\cA^\bullet_M(\bbW_M),\nabla)$ is a sheaf of differential graded 
algebras over the $C^\infty$ de Rham complex $(\cA^\bullet_M,d)$
(in particular $H^\bullet(\cA^\bullet_M(\bbW_M),\nabla)$ is a sheaf of
graded algebras over the constant sheaf $\bbC_M[[\hbar]]$);

\item there is a filtered quasiisomorphism
\[
(\bbA^\hbar_M,F_\bullet) @>>> ((\cA^\bullet_M(\bbW_M),\nabla),F_\bullet)
\]
of differential graded algebras over $\bbC_M[[\hbar]]$ (in particular
$H^p(F_q\cA^\bullet_M(\bbW_M),\nabla)=0$ for $p\neq 0$);
\end{enumerate} 
therefore $\nabla$ is determined by its component
\[
\nabla : \cA^0_M(\bbW_M) @>>> \cA^1_M(\bbW_M)
\]
which has all the properties of a flat connection on $\bbW_M$.

Let $H = Sp(\dim M)$ and let $P @>{\pi}>> M$ denote the $H$-principal bundle
of symplectic frames in $TM$, identify $TM$ with the vector bundle assiciated
to the standard representation of $H$.

Recall that, for an $\frak{h}$-module $V$, the subcomplex
$\left[\pi_*\cA_P^\bullet\otimes V\right]^{basic}\subset\pi_*\cA_P^\bullet$
is defined by the pull-back diagram
\[
\begin{CD}
\left[\pi_*\cA_P^\bullet\otimes V\right]^{basic} @>>>
C^\bullet(\pi_*\cT_P,\frak{h};V) \\
@VVV  @VVV \\
\pi_*\cA_P^\bullet\otimes V @>>> C^\bullet(\pi_*\cT_P;V)
\end{CD}
\]
where $\cT_P$ denotes the sheaf of Lie algebras of $C^\infty$ vector fields
on $P$.

Let $W$ denote the Weyl algebra
of the standard representation of $H$. Then $\bbW_M$ is identified with
the sheaf of sections of the associated bundle $P\times_HW$ and
pull-back by $\pi$
\[
\cA_M^\bullet(\bbW_M) @>>> \pi_*\cA_P^\bullet\widehat\otimes W
\]
identifies $\bbW_M$-valued forms on $M$ with the subcomplex of basic
$W$-valued forms on $P$. The flat connection $\nabla$ gives rise to
a basic $\frak{g}(=Der(W))$-valued 1-form
$A\in H^0(P;\cA_P^1\widehat\otimes\frak{g})$ which satisfies the
Maurer-Cartan equation $dA+\frac12\lbrack A,A\rbrack = 0$ (so that
$(d + A)^2=0$). Then, pull back by $\pi$ induces the isomorphism of
filtered complexes
\[
(\cA_M^\bullet(\bbW_M),\nabla) @>>>
\left(\left\lbrack\pi_*\cA_P^\bullet\widehat\otimes W\right\rbrack^{basic},
d+A\right)\ .
\]

Given $A$ as above and a (filtered) topological $\frak{g}$-module $L$ such
that the action of $\frak{h}\subset\frak{g}$ integrates to an action of $H$.
Set
\[
(\cA_M^\bullet(L),\nabla)\overset{def}{=}
\left(\left\lbrack\pi_*\cA_P^\bullet\widehat\otimes L\right\rbrack^{basic},
d+A\right)\ .
\]

Note that the association $L\mapsto (\cA_M^\bullet(L),\nabla)$ is functorial
in $L$. In particular it extends to complexes of $\frak{g}$-modules.

Taking $L=\bbC$, the trivial $\frak{g}$-module, we recover
$(\cA_M^\bullet,d)$. For any complex $(L^\bullet,d_L)$ of $\frak{g}$-modules
as above the complex $(\cA_M^\bullet(L^\bullet),\nabla + d_L)$ has a natural
structure of a differential graded module over $(\cA_M^\bullet,d)$.

The Gel'fand-Fuchs map
\[
GF : C^\bullet(\frak{g},\frak{h};L) @>>> \cA_M^\bullet(L)
\]
(the sourse understood to be the constant sheaf) is defined by the
formula
\[
GF(c)(X_1,\ldots,X_p) = c(A(X_1),\ldots,A(X_p))\ ,
\]
where $c\in C^p(\frak{g},\frak{h};L)$ and $X_1,\ldots,X_p$ are locally
defined vector fields. It is easy to verify that $GF$ takes values in
basic forms and is a map of complexes. Note also, that $GF$ is natural
in $L$. In particular the definition above has an obvious extension to
complexes of $\frak{g}$-modules.

We now proceed to apply the above constructions to particular examples
of complexes $L^\bullet$ of $\frak{g}$-modules. In all examples below
the $\frak{g}$-modules which appear have the following additional
property which is easy to verify, namely,
\[
\text{$H^p(\cA_M^\bullet(L),\nabla) = 0$ for $p\neq 0$}\ .
\]
If $(L^\bullet,d_L)$ is a complex of $\frak{g}$-modules with the above
property, then the inclusion
\begin{equation}\label{map:incl-ker}
(\ker(\nabla),d_L\vert_{\ker(\nabla)})\hookrightarrow
(\cA_M^\bullet(L^\bullet),\nabla + d_L)
\end{equation}
is a quasiisomorphism. In such a case the map \eqref{map:incl-ker} induces
an isomoprhism
\[
\R\Gamma(M;\ker(\nabla)) @>>> \Gamma(M;\cA_M^\bullet(L^\bullet))
\]
in the derived category of complexes since the sheaves
$\cA_M^\bullet(L^\bullet)$ are soft and the Gel'fand-Fuchs map induces
the (natural in $L^\bullet$) morphism
\begin{equation}\label{map:GF}
GF : C^\bullet(\frak{g},\frak{h};L) @>>> \R\Gamma(M;\ker(\nabla))
\end{equation}
in the derived category.

For $L = W$ (respectively $\frak{g}$, $CC^{per}_\bullet(W)$, $\widehat\cO$,
$\bbC$, $\widehat\Omega^\bullet$, $CC^{per}_\bullet(\widehat\cO)$),
$\ker(\nabla) = \bbA^\hbar_M$ (respectively $Der(\bbA^\hbar_M)$,
$CC^{per}_\bullet(\bbA^\hbar_M)$, $\cO_M$, $\bbC_M$, $\Omega^\bullet_M$,
$CC^{per}_\bullet(\cO_M)$). We leave it to the
reader to identify $\ker(\nabla)$ for other (complexes of) $\frak{g}$-modules
which appear in Section \ref{section:formalRR} by analogy with the above
examples.

The relationship between the Lie algebra cocycles defined in Section
\ref{section:formalRR} and morphism defined in Section
\ref{section:cx-man} established by the Gel'fand-Fuchs map
\eqref{map:GF} is as follows.

\subsubsection{The trace density}
The image of $\tilde\mu^\hbar$ (defined in Lemma \ref{lemma:formalTR})
under (the map on cohomology in degree zero induced by)
\begin{multline*}
GF : C^\bullet(\frak{g},\frak{h};\Hom^\bullet
(CC^{per}_\bullet(W)[\hbar^{-1}],\prod_{p\in\bbZ}\widehat\Omega^\bullet
[\hbar^{-1},\hbar]][2d-2p]) @>>> \\
\R\Gamma(M;\shHom^\bullet(CC^{per}_\bullet(\bbA^\hbar_M)[\hbar^{-1}],
\prod_{p\in\bbZ}\Omega^\bullet_M[\hbar^{-1},\hbar]][2d-2p]))
\end{multline*}
is the morphism $\tilde\mu^\hbar_{\bbA}$ defined in \eqref{map:TR}. 

Similarly, the image of $\tilde\mu^\hbar\smile\iota$ under $GF$ is the
morphism $\tilde\mu^\hbar_{\bbA}\circ\iota$.

\subsubsection{The symbol and the Hochschild-Kostant-Rosenberg map}
The image of $\sigma$ (defined in \ref{sssection:HKR})
under (the map on cohomology in degree zero induced by)
\[
GF : C^\bullet(\frak{g},\frak{h};
\Hom^\bullet(CC^{per}_\bullet(W), CC^{per}_\bullet(\widehat\cO))
@>>>
\R\Gamma(M;\shHom^\bullet(CC^{per}_\bullet(\bbA^\hbar_M),
CC^{per}\bullet(\cO_M))
\]
is the morphism
\[
\sigma : CC^{per}_\bullet(\bbA^\hbar_M) @>>> CC^{per}_\bullet(\cO_M)\ .
\]
The image of $\tilde\mu$ (defined in \ref{sssection:HKR})
under (the map on cohomology in degree zero induced by)
\begin{multline*}
GF : C^\bullet(\frak{g},\frak{h};\Hom^\bullet
(CC^{per}_\bullet(\widehat\cO),\prod_{p\in\bbZ}\widehat\Omega^\bullet[-2p]))
@>>> \\
\R\Gamma(M;\shHom^\bullet(CC^{per}_\bullet(\cO_M),
\prod_{p\in\bbZ}\Omega^\bullet_M[-2p]))
\end{multline*}
is the morphism $\tilde\mu_{\cO}$ defined in \eqref{map:HKR}. 

\subsubsection{The characteristic class of the deformation}
The image of the cocycle $\theta\in C^2(\frak{g},\frak{h};
\frac{1}{\sqrt{-1}\hbar}\bbC[[\hbar]])$ (defined in \ref{sssection:charcl})
under the map
\[
GF: C^\bullet(\frak{g},\frak{h};\frac{1}{\sqrt{-1}\hbar}\bbC[[\hbar]]) @>>>
\frac{1}{\sqrt{-1}\hbar}\cA^\bullet_M[[\hbar]]
\]
is the characteristic class $\theta$ of the deformation quantization
$\bbA^\hbar_M$ defined in \cite{F} and \cite{D}.

\subsubsection{The $\widehat A$-class}
The composition
\[
\widehat\Symm^\bullet(\frak{h})^{\frak{h}} @>{CW}>>
C^\bullet(\frak{g},\frak{h};\bbC) @>{GF}>> \cA_M^\bullet
\]
is easily seen to be the usual Chern-Weil homomorphism. In particular
we have
\[
GF(\widehat A) = \widehat A(TM)\ ,
\]
where $\widehat A$ is defined in \ref{sssection:formalA}.

Combining the above facts we obtain the following proposition.

\begin{prop}
The image of the cocycle
$\tilde\mu^\hbar\smile\iota - \tilde\mu\smile\widehat
A\smile e^\theta\smile\sigma$ (see Theorem
\ref{thm:formalRR}) is the morphism
$\iota\circ\tilde\mu^\hbar_{\bbA}-(\tilde\mu_{\cO}\widehat
A\smile e^\theta)\circ\sigma$ (see Theorem
\ref{thm:main}).
\end{prop}

\begin{cor}\label{cor:reduction}
Theorem \ref{thm:formalRR} implies Theorem \ref{thm:main}.
\end{cor}


\begin{thebibliography}{ABCDEF}
\bibitem[BFFLS]{BFFLS} F.~Bayen, M.~Flato, C.~Fronsdal, A.~Lichneroviz
and D.~Sternheimer, {\em Deformation theory and quantizaton {\rm I} and
{\rm II}}, Annals of Physics {\bf 111} (1978), 61--151.
\bibitem[B]{B} J.-E.~Bjork, {\em Analytic $\cD$-modules and Applications},
Kluwer Academic Publishers, 1993.
\bibitem[Bry]{Bry} J.-L.~Brylinski, {\em Some examples of Hochschild and
cyclic homology}, Utrecht 1986, A.~Cohen et al. eds., Lecture Notes in Math.
vol. 1271, Springer-Verlag 1987, 33--72.
\bibitem[D]{D} P.~Deligne, {\em deformattions de l'algebre des fonctions
d'une variete symplectique, Comparison entre Fedosov et De Wilde, Lecomte},
Selecta Math., New series {\bf 1} (1995), 667--698.
\bibitem[F]{F} B.V.~Fedosov, {\em Deformation quantizaton and index theory},
Mathematical Topics 9, Academie Verlag, Berlin, 1996. 
\bibitem[FT1]{FT1} B.~Feigin and B.~Tsygan, {\em Cohomology of Lie algebras
of generalized Jacobi matices}, Funct. Anal. and Appl. {\bf 17} (1983).
\bibitem[FT2]{FT2} B.~Feigin and B.~Tsygan, {\em Lie algebra homology and
the Riemann-Roch theorem}, Proccedings of the winter school of 2nd winter
school at Srni, Rend. Math. Palermo (1989).
\bibitem[Ill]{Ill} L.~Illusie, {\em G\'en\'eralit\'es sur les conditions de
finitude dans les cat\'egories deriv\'ees}, in SGA6, Lecture Notes in
Math. vol. 225, Springer-Verlag, 1971, 78--159.
\bibitem[KS]{KS} M.~Kashiwara and P.~Schapira, {\em Sheaves on manifolds},
Grundlehren der mathimatischen Wissenschaften, no. 292, Springer-Verlag, 1990.
\bibitem[McC]{McC} R.~McCarthy, {\em The cyclic homology of an exact category},
Journal of Pure and Applied Algebra {\bf 93} (1994), 251--296.
\bibitem[M]{M} S.A.~Mitchell, {\em Hypercohomology spectra and Thomason's
Descent theorem}, preprint, 1996.
\bibitem[NT1]{NT1} R.~Nest and B.~Tsygan, {\em Algebraic index theorem},
Comm. Math. Phys. {\bf 172} (1995), 223--262.
\bibitem[NT2]{NT2} R.~Nest and B.~Tsygan, {\em On the cohomology ring of
an algebra}, preprint, 1996,
\bibitem[NT3]{NT3} R.~Nest and B.~Tsygan, {\em Deformation quantization of
symplectic Lie algebroids}, preprint, 1997.
\bibitem[SS]{SS} P.~Schapira and J.-P.~Schneiders, {\em Index theorm for
elliptic pairs, Elliptic pairs} II, Asterisque vol. 224, 1994.
\bibitem[T]{T} R.~Thomason, {\em Algebraic K-theory and \'etale cohomology},
Ann. Scient. Ecole Norm. Sup. {\bf 13} (1985), 393--444.
\bibitem[TT]{TT} R.~Thomason and T.~Trobaugh, {\em Higher algebraic K-theory of
schemes and of derived categories}, The Grothendieck Festschrift, Volume
III, Progress in Mathematics Vol. 88, Birk\"auser 1990, 247-436.
\bibitem[W]{W} F.~Waldhausen, {\em Algebraic K-theory of spaces}, Lecture notes
in Math. vol. 854, Springer-Verlag, 1981, 494-517.
\end{thebibliography}
\end{document}